# Vector Determinant in the ARE Framework: From Scalar to Vector-Valued

**Author:** Ramón Moya

**Affiliation:** School of Mathematics, Autonomous University of Santo Domingo (UASD), Dominican Republic

**E-mail:** rmoya07@uasd.edu.do

**ORCID:** 0009-0001-1601-4699

**Abstract**

**The ARE (Action, Rectification, and Structure)** method is presented as a framework for reorganizing the terms of the Leibniz expansion before taking the scalar sum that yields the determinant. The starting point is the right action of the cyclic group $C_n$ on $S_n$, which partitions the permutations into orbits of size $n$. By fixing a canonical representative per orbit, natural orbital classes are obtained, and from these, a finite spectral decomposition is derived using a Discrete Fourier Transform.

The central object of the work is built upon this organization: the vector determinant

$D_\phi(A) = (G_0(A), G_1(A), \dots, G_{n-1}(A)),$

where each mode $G_k(A)$ is the Discrete Fourier Transform of the orbital sums $\Lambda_r(A)$. The classical determinant appears exactly as the fundamental mode $G_0(A) = \det(A)$. For $k \neq 0$, the modes $G_k$ are generally not alternating, so they should be understood as multilinear spectral observables associated with the orbital structure rather than as additional classical determinants.

The manuscript rigorously establishes the fundamental properties of this formalism: exact recovery of the determinant by means of a linear readout functional, multilinearity in all modes, Hermitian symmetry for real matrices, elementary bounds, vector Jacobi and Laplace formulas, an orbital Parseval identity, a vector Hadamard inequality, and finite trigonometric interpolation.

This construction admits a continuous extension to all real frequencies via the finite trigonometric interpolation

$$G(z;A) = \sum_{r=0}^{n-1} \Lambda_r(A)\, e^{2\pi i r z/n},$$

which extends the discrete modes to continuous frequencies and leads to the associated orbital polynomial.

The approach is structural and conceptual: it does not reduce the asymptotic factorial complexity of the general calculus of dense determinants, but it does reveal an orbital and spectral layer that the scalar $\det(A)$ does not distinguish. In this version, the emphasis is on the spectral core of the ARE program. Additional issues of the broader geometric framework, including orbital cancellations due to symmetry and dihedral-geometric developments, should be understood as part of the general ARE program and not as central results established in this article.

The work includes fully computed low-dimensional examples, operational formulations of the formalism, and computational reproducibility material deposited on Zenodo (DOI: 10.5281/zenodo.17423738).

# 1. Introduction

## 1.1. Motivation and context

The determinant of a matrix $n \times n$ is classically defined as the sum of $n!$ Leibniz terms, each a signed product of $n$ entries. For $n = 3$, Sarrus's rule visually organizes the 6 terms diagonally; for $n \geq 4$, any attempt to directly extend that diagram encounters an unavoidable combinatorial barrier: the number of terms grows factorially, while the geometric width of a Sarrus-type scheme remains rigid. The question, therefore, is not how to literally extend Sarrus's drawing, but what underlying algebraic structure explains its effectiveness in $3 \times 3$ and how to generalize that structure to all $n$.

The starting point of the ARE Method is precisely that observation: Leibniz's terms do not appear in isolation, but are organized in orbits under the action of the cycle

$\rho = (1\ 2\ \cdots\ n).$

Under the composition on the right, this action partitions $S_n$ into $(n-1)!$ orbits of cardinality $n$. Each orbit admits a canonical representative $\sigma^*$ with $\sigma^*(1) = 1$, which allows for the clear separation of two combinatorial levels: the free permutation of the columns $2, \dots, n$, and the cyclic rotation that completes the orbit. This observation naturally organizes the Leibniz expansion before the final summation and constitutes the first pillar of the ARE program.

The second pillar is geometric. Canonical rectification transforms each orbit into a uniform visual pattern, where cyclic rotation becomes a rigid translation on the grid. This rectification does not reduce combinatorial complexity, but it does reveal a structural regularity that remains hidden in the usual scalar formula. The third pillar is the dihedral structure, which pairs orbits through involution $\Phi$ and explains robust cancellations in matrices with discrete symmetries. Thus, ARE should not be understood as an asymptotically faster algorithm, but rather as an algebraic-geometric rearrangement of the determinant that isolates patterns of symmetry, cancellation, and phase before scalar collapse.

The crucial connection to the present article emerges when we observe that the cyclic rotation of an orbit is algebraically analogous to the cyclic translation that diagonalizes the **Discrete Fourier Transform**. This coincidence makes it natural to introduce, for each matrix $A$, the orbital sums

$$\Lambda_r(A) = \sum_{\phi(\sigma)=r} m_\sigma(A)$$

and its Fourier modes

$$G_k(A) = \sum_{r=0}^{n-1} \Lambda_r(A)\, \omega^{kr}$$

Thus, the classical determinant is no longer the only observable, and it appears as the fundamental mode:

$$G_0(A) = \det(A).$$

The resulting object,

$$D_\phi(A) = (G_0(A), G_1(A), \dots, G_{n-1}(A)).$$

is the vector determinant and constitutes the main focus of the manuscript.

It is important to emphasize from the outset an essential conceptual distinction. The vector $D_\varphi(A)$ enriches the determinant data of A, but its components do not all play the same role: the mode $G_0(A)$ coincides with $det(A)$, while the modes $G_k(A)$, $k \neq 0$, are multilinear observables associated with the orbital structure and, in general, are not alternating. Consequently, the term “vector determinant” should be understood as a spectral extension of the determinant data, not as a family of independent classical determinants. This clarification avoids confusing the alternating scalar level of the determinant with the additional modal layer introduced by the ARE formalism.

From this perspective, the problem changes in nature. It is no longer simply a matter of calculating the determinant $\det(A)$, but of understanding what additional information about the orbital distribution of the monomials is encoded in the vector determinant $D_{\phi}(A)$. Two matrices with the same determinant can have radically different orbital configurations; the classical scalar does not distinguish these cases, while the vector determinant $D_{\phi}(A)$ can. This observation places the vector determinant at the intersection of multilinear algebra, finite harmonic analysis, permutation combinatorics, and discrete geometry.

The aim of this work is to rigorously establish the spectral layer of the ARE program. In particular, we develop the framework

$$M\phi(A) \longrightarrow \{\Lambda r(A)\} \leftrightarrow D\phi(A) \longrightarrow G(z;A)$$

where $M_{\phi}(A)$ is the vector of phased monomials, $\{\Lambda_r(A)\}$ are the orbital sums, $D_{\phi}(A)$ is its Fourier transform, and $G(z;A)$ is the associated finite trigonometric interpolation. The result is a formalism that maintains full compatibility with the classical determinant, but extends its interpretation to a theory of modes, energies, correlations, and spectral zeros.

This article should therefore be read as the rigorous formulation of the spectral level of the ARE program. The purely orbital and geometric part motivates the construction, whereas the vector and analytic part develops it rigorously. Complementary works in the series delve deeper, on the one hand, into the ARE Method as an autonomous combinatorial-geometric framework and, on the other hand, into possible tensor and hyperdeterminant extensions.

For further developments of the broader geometric-orbital framework, see [6].

Before describing the organization of the manuscript, it is useful to state explicitly the scope and level of originality of the present contribution. The framework developed here combines multilinear algebra, permutation group actions, and discrete harmonic analysis, and therefore touches several existing mathematical traditions. To avoid ambiguity regarding what is genuinely introduced in this work versus what follows naturally from the formalism, we summarize below the principal contributions, limitations, and intended scope of the theory.

### 1.2. Novelty and Scope of the Present Work

This article introduces a vector-valued multilinear formalism associated with matrices, called the vector determinant,

$$D_\varphi(A) = (G_0(A), \dots, G_{n-1}(A)),$$

constructed from the discrete Fourier transform of orbital sums associated with the Leibniz expansion.

The originality of the work lies primarily in the following aspects:

1. The systematic definition of the vector determinant as an algebraic object attached to arbitrary matrices.
2. The construction of orbital modes $G_k(A)$ indexed by the characters of the cyclic group $C_n$.
3. The proof that matrices with equal determinant may possess distinct determinant vectors.
4. The introduction of an orbital operator whose spectrum is exactly diagonalized by the discrete Fourier transform.
5. The degree-incompatibility result showing that the modal quantities $G_k(A)$ cannot generally be identified with eigenvalues of circulant matrices.

Several additional results, including multilinearity, vectorial Jacobi identities, orbital Parseval identities, and modal inequalities, are developed as structural consequences of this framework.

The present work does not claim to solve inverse synthesis, injectivity, or probabilistic concentration problems, which remain open directions.

### 1.3. Relation with Classical Determinants and Harmonic Analysis

The present construction differs fundamentally from several classical frameworks.

The determinant of Dedekind–Frobenius associates a polynomial invariant to group-indexed matrices, particularly for abelian groups. In contrast, the determinant vector introduced here applies to arbitrary matrices and is derived from the orbital organization of Leibniz monomials.

Likewise, although the discrete Fourier transform diagonalizes circulant matrices, the modal quantities $G_k(A)$ do not coincide with circulant eigenvalues in general.

The vector determinant therefore occupies an intermediate position between multilinear algebra, permutation group actions, and harmonic analysis on finite cyclic structures.

### 1.4. Article Structure

Section 2 introduces the space V, the phase function φ, the phased monomial vector $M_{\varphi}(A)$, and the readout functional $L_{\varphi}$, culminating in Corollary 2.6 (Uniqueness of the readout functional). Section 3 defines orbital sums $\Lambda_r(A)$, modes $G_k(A)$, and the determinant vector $D_{\varphi}(A)$, and establishes their basic properties. Section 4 presents fully computed examples in low dimensions. Section 5 introduces intrinsic vector minors, modal cofactors, and the vector Jacobi and Laplace formulas. Section 6 studies Fourier compatibility and the circulant case, including the orbital multiplication operator and the structural limit imposed by degree incompatibility, and shows that there is no direct identification between the determinant vector $D_{\phi}(A)$ and the spectrum of circulant matrices associated with it $A$, due to polynomial-degree incompatibility. Section 7 clarifies two conceptual corrections: the distinction between the diagonal product and the trace, and the merely motivational scope of the connections with physics. Section 8 contains probabilistic results for real iid Gaussian matrices and a table of open problems. Section 9 develops vector operations on $D_{\phi}(A)$ orbital Parseval, autocorrelation, modal energy, modal cofactor matrix, entropy, and the vector Hadamard inequality. Section 10 studies finite trigonometric interpolation $G(z; A)$, the orbital polynomial $P_A(q)$, and associated spectral zeros. Section 11 formalizes orbital phase canonicity, covariance under global affine gauge, and genuine dependence under orbital gauge. Finally, the conclusions precisely delimit established results, current limitations of the formalism, and open directions for research.

## Main notation

| Symbol | Meaning |
|---|---|
| $M_n(\mathbb{C})$ | Matrices with complex entries |
| $S_n$ | Symmetric group of permutations of $\{1, \dots, n\}$ |
| $\rho$ | Standard cycle generator $(1\ 2\ \cdots\ n)$ |
| $\omega$ | Primitive $n$-th root of unity (positive convention) |
| $\varepsilon(\sigma)$ | Sign of the permutation $\sigma$ |
| $m_\sigma(A)$ | Leibniz term associated with $\sigma$ |
| $\phi(\sigma)$ | Orbital phase function |
| $\Lambda_r(A)$ | Orbital sum at phase $r$ |
| $M_\phi(A)$ | Vector of phased monomials |
| $L_\phi$ | Readout functional |
| $G_k(A)$ | Orbital mode at frequency $k$ |
| $D_\phi(A)$ | Determinant vector |
| $G(z; A)$ | Interpolated function |
| $P_A(q)$ | Orbital polynomial |
| $\text{per}\,(A)$ | Permanent of $A$ |
| $C_{ij}^{(k)}(A)$ | $k$-modal cofactor |
| $C^{(k)}(A)$ | $k$-modal cofactor matrix |
| $H_A$ | Orbital multiplication operator |

**Theorem A (Synthesis).**

Given the canonical orbital phase $\phi$, the determinant vector $D_\phi(A) = (G_0(A), \dots, G_{n-1}(A)) \in \mathbb{C}^n$ satisfies, for all $A \in M_n(\mathbb{C})$:

(i) $G_0(A) = \det\,(A)$.

(ii) Each mode $G_k(A)$ is multilinear in rows and columns of $A$.

(iii) For $k \neq 0$, $G_k(A)$ it is not alternating in general.

(iv) If $A \in M_n(\mathbb{R})$, then $G_{n-k}(A) = \overline{G_k(A)}$ for all $k$.

(v) $\frac{1}{n} \parallel D_\phi(A) \parallel_2^2 = \sum_{r=0}^{n-1} \mid \Lambda_r(A) \mid^2$ (Orbital Parseval).

(vi) Jacobi and Laplace-type formulas are verified mode by mode via modal cofactors.

(vii) The classical spectrum is $A$ recovered not from the crude modes $G_k(A)$ of total order, but from the hierarchy of fundamental modes $G_0^{(r)}(A[I,I])$ of the principal minor orders $r = 1,\dots,n$.

The full proofs are given in §§3, 5, 9 and 11.

**1.5. Precise scope of the article**

**Proven results:** This work establishes the definition and fundamental properties of the determinant vector $D_\phi(A)$: exact recovery of $\det(A)$, multilinearity, non-alternation for $k \neq 0$, Hermitian symmetry for real matrices, orbital Parseval, Jacobi and Laplace modal formulas, and the vector Hadamard inequality.

**Corrective structural result:** It is proven that crude modes $G_k(A)$ cannot be identified with the eigenvalues of circulant matrices due to algebraic incompatibility of degrees, and that the correct link to the classical spectrum passes through the fundamental modes of the hierarchy of principal minors. In particular, there is no naturally constructed circulant operator from $A$ whose list of eigenvalues is $(G_0(A),\dots,G_{n-1}(A))$: such an identification is impossible for $n \geq 2$ due to irreducible incompatibility of polynomial degrees (Theorem 11.12). The correct link between the vector determinant and the classical spectrum is established exclusively through the fundamental modes $G_0^{(r)}(A[I,I])$of the principal minors, as formalized in Section 5.

**Result of distinction:** A pair with $\det(A) = \det(B)$ but $D_\phi(A) \neq D_\phi(B)$ is explicitly exhibited $A, B$, illustrating that the determinant vector carries orbital information not accessible to the classical scalar.

**Open questions:** The following remain as open problems: the inverse synthesis (surjectivity of the map $A \mapsto \{\Lambda_r(A)\}$), the characterization of the zeros of the orbital polynomial $P_A$ on the unit circle, the Gaussian concentration of $H(A) \to \log\, n$ in probability, and the possible tensor extensions of the formalism.

**1.6. Dependence on the orbital phase**

Throughout the work, a canonical orbital phase is established $\phi: S_n \to \mathbb{Z}_n$ (Definition 2.2), constructed from representatives with $\sigma^*(1) = 1$ and $\phi(\sigma^*) = 0$. It is important to distinguish which objects are intrinsic and which depend on the choice of $\phi$.

**Independent objects of $\phi$:** the classical determinant $\det(A)$; the family of classical minors $\det(A[I,J])$; the quotient $R_\infty(A) = |\det(A)|/\mathrm{per}(|A|)$ (Definition 3.10).

**Objects that depend on the choice of $\phi$:** the phased monomial vector $M_\phi(A)$ (Definition 2.3); the orbital sums $\Lambda_r(A)$ (Definition 3.1); the determinant vector $D_\phi(A)$ (Definition 3.3); the modal variance $V_2(A)$ and the spectral entropy $H(A)$ (Definitions 3.10 and 9.10).

All results in this article are stated for the canonical phase $\phi$. Changing the phase corresponds to a change in orbital gauge that conserves $\det(A)$ but redistributes the energy between the modes $G_k(A)$; the complete structure of this phenomenon is formalized in §11.

## 2. Space of Phased Monomials and Recovery of the Determinant

### 2.1. Fundamental Objects

**Definition 2.1 (Space of orbital functions).**

Let

$$V := \mathbb{C}^{S_n}$$

be the complex vector space of all functions on the symmetric group $S_n$. We endow $V$ with its canonical basis $\{\delta_\sigma\}_{\sigma\in S_n}$, where

$$\delta_\sigma(\tau) = \begin{cases} 1, & \text{if } \tau = \sigma, \\ 0, & \text{otherwise.} \end{cases}$$

Then

$$\dim_{\mathbb{C}} V = n!, \dim_{\mathbb{R}} V = 2\,n!$$

For the general theory of function spaces on finite groups, see [5].

**Definition 2.2 (Orbital phase function)**

A phase function is a map $\phi: S_n \to \mathbb{Z}_n$ such that, for each orbit $O(\sigma^*) = \{\sigma^* \circ \rho^r : r \in \mathbb{Z}_n\}$, there exists a unique representative $\sigma^*$ with

$$\sigma^*(1) = 1, \phi(\sigma^*) = 0, \phi(\sigma^* \circ \rho^r) = r \quad (r \in \mathbb{Z}_n).$$

**Algorithmic construction:**

Given $\sigma \in S_n$, its canonical representative is

$$\sigma^* = \sigma \circ \rho^{r^*}, r^* \equiv \sigma^{-1}(1) - 1 (\mathrm{mod}\ n),$$

verifying that $\sigma^*(1) = (\sigma \circ \rho^{r^*})(1) = \sigma(r^* + 1) = \sigma(\sigma^{-1}(1)) = 1$. The position of $\sigma$ in its orbit is $\phi(\sigma) = r^*$.

**Comment.** Definition 2.2 encodes in $\phi$ the orbital structure of the action of the cyclic subgroup $C_n = \langle\rho\rangle$ on $S_n$ by right-hand composition: $(\sigma, r) \mapsto \sigma \circ \rho^r$. Each orbit contains exactly $n$elements and contributes a single element to each class $\phi^{-1}(r)$, which will be essential for defining the orbital sums $\Lambda_r(A)$ and their Fourier modes $G_k(A)$ (Definitions 3.1 and 3.2).

**Definition 2.3 (Phased monomial vector)**

Once fixed $\phi$, we define

$$M_\phi(A) \coloneqq \sum_{\sigma \in S_n} m_\sigma(A) \cdot \omega^{\phi(\sigma)} \cdot \delta_\sigma \ \in \ V$$

Here $\delta_\sigma$ denotes the canonical basis of $V = \mathbb{C}^{S_n}$, given by $\delta_\sigma(\tau) = 1$ if $\tau = \sigma$ and 0 otherwise; therefore, $M_\phi(A)$ it is simply the vector whose coordinates are $m_\sigma(A)\omega^{\phi(\sigma)}$.

**Definition 2.4 (Readout functional)**

The linear functional associated with $\phi$ is

$$L_\phi : V \to \mathbb{C}, L_\phi(f) := \sum_{\sigma \in S_n} \omega^{-\phi(\sigma)} f(\sigma).$$

**2.2. Exact recovery of the determinant**

**Theorem 2.5 (Recovery)**

For all $A \in M_n(\mathbb{C})$,

$\det(A) = L_\phi(M_\phi(A))$.

Proof:

By definition of functional $L_\phi$ (Definition 2.4):

$$L_\phi(M_\phi(A)) = \sum_{\sigma \in S_n} \omega^{-\phi(\sigma)} \ M_\phi(A)(\sigma).$$

The coordinate $\sigma$ of the vector $M_\phi(A)$ is, by Definition 2.3:

$M_\phi(A)(\sigma) = m_\sigma(A)\, \omega^{\phi(\sigma)}$.

Substituting and cancelling the phase factors:

$$L_\phi\left(M_\phi(A)\right) = \sum_{\sigma \in S_n} \omega^{-\phi(\sigma)} \cdot m_\sigma(A) \cdot \omega^{\phi(\sigma)} = \sum_{\sigma \in S_n} m_\sigma(A) \cdot \underbrace{\omega^{\phi(\sigma)-\phi(\sigma)}}_{=1} = \sum_{\sigma \in S_n} m_\sigma(A)$$

$= \det(A)$

**Observation 2.5.1 (Purely algebraic character of the formalism).** In the definition of $M_\phi(A)$ and in the functional $L_\phi$, no notion of analytic convergence intervenes: all sums are finite (over the finite set $S_n$) and the formalism is entirely algebraic.

**2.3. Uniqueness of the Readout Functional**

**Corollary 2.6 (Uniqueness of the Readout Functional).**

Let $F: V \to \mathbb{C}$ be a linear functional such that

$F(M_\phi(A)) = \det(A)$ for all $A \in M_n(\mathbb{C})$.

Then $F = L_\phi$.

**Proof.** Write

$$F(f) = \sum_{\sigma \in S_n} c_\sigma \ f(\sigma), c_\sigma \in \mathbb{C}.$$

**Step 1.** Fix $\sigma_0 \in S_n$, and consider the permutation matrix $P_{\sigma_0}$ defined by

$$\left(P_{\sigma_0}\right)_{ij} = \begin{cases} 1, & \text{if } j = \sigma_0(i), \\ 0, & \text{otherwise.} \end{cases}$$

For this matrix, the only nonzero Leibniz term is the one corresponding to $\sigma_0$:

$m_{\sigma_0}(P_{\sigma_0}) = \varepsilon(\sigma_0), m_\sigma(P_{\sigma_0}) = 0$ for $\sigma \neq \sigma_0$.

Therefore,

$$M_\phi(P_{\sigma_0}) = \varepsilon(\sigma_0)\, \omega^{\phi(\sigma_0)}\, \delta_{\sigma_0}.$$

**Step 2.** Applying $F$ and using the hypothesis, we obtain

$$\varepsilon(\sigma_0) = \det(P_{\sigma_0}) = F(M_\phi(P_{\sigma_0})) = c_{\sigma_0}\, \varepsilon(\sigma_0)\, \omega^{\phi(\sigma_0)}.$$

Since $\varepsilon(\sigma_0) \neq 0$, it follows that

$$c_{\sigma_0} = \omega^{-\phi(\sigma_0)}.$$

**Step 3.** Since $\sigma_0 \in S_n$ was arbitrary, we obtain

$c_\sigma = \omega^{-\phi(\sigma)}$ for all $\sigma \in S_n$.

Because $\{\delta_\sigma\}_{\sigma \in S_n}$ is a basis of $V$, the functional $F$ is uniquely determined, and hence

$$F = L_\phi.$$

In particular, the identity $F = L_\phi$ holds not only on the subfamily

$\{M_\phi(A): A \in M_n(\mathbb{C})\}$,

but on the whole space $V$. That is, any linear functional that recovers $\det(A)$ from the phased monomial vector $M_\phi(A)$ must coincide with $L_\phi$ on every basis vector.

**Conceptual comment.** This corollary sharply identifies the novelty of the formalism: the determinant cannot be recovered from $M_\phi(A)$ by any linear procedure other than undoing the monomial phases term by term. The additional information is therefore encoded in the object $M_\phi(A)$ itself and in the new invariants that can be constructed from it.

This is a rigidity phenomenon for linear functionals.

## 3. Orbital Sums, Modes and Basic Properties

### 3.1. The hierarchy of objects

**Definition 3.1 (Orbital sums)**

For $r = 0,1, \dots, n-1$:

$$\Lambda_r(A) := \sum_{\substack{\sigma \in S_n \\ \phi(\sigma) = r}} m_\sigma(A). \qquad (3)$$

Each $\Lambda_r(A)$ sums all Leibniz terms whose orbital position is $r$. For each $r$ there are exactly $(n-1)!$ permutations with $\phi(\sigma) = r$ (one for each orbit). If $A \in M_n(\mathbb{R})$, then $\Lambda_r(A) \in \mathbb{R}$ for all $r$.

**Lemma 3.1bis (Exact cardinality of each orbital class).**

For each

$r \in \mathbb{Z}/n\mathbb{Z}$,

the preimage

$\phi^{-1}(r) = \{\sigma \in S_n : \phi(\sigma) = r\}$

contains exactly $(n-1)!$ permutations.

**Proof.** Let

$r \in \mathbb{Z}/n\mathbb{Z}$

be fixed. By Definition 2.2, a permutation $\sigma$ satisfies

$\phi(\sigma) = r$

if and only if it occupies position $r$within its orbit, that is, if and only if

$\sigma = \sigma^* \circ \rho^r$

for some canonical representative $\sigma^*$. Since the canonical representatives are in bijection with the orbits, and there are exactly $(n-1)!$ orbits (Proposition 6.1, proved in §6), the class

$\phi^{-1}(r)$

has exactly $(n-1)!$ elements.

**Direct verification by counting.** The total is $\sum_{r=0}^{n-1} |\, \phi^{-1}(r) \,| = n \cdot (n-1)! = n! = |\, S_n \,|$, which is consistent. For $n = 3$: $(3-1)! = 2$ permutations by class; for $n = 4$: $(4-1)! = 6$ permutations by class. This lemma guarantees that orbital sums $\Lambda_r(A)$ always add up to the same number $(n-1)!$ of monomials, giving them precise combinatorial symmetry.

**Corollary 3.1.ter.** The identity $\sum_{r=0}^{n-1} \Lambda_r(A) = \det(A)$ is obtained by summing over all classes:

$$\sum_{r=0}^{n-1} \Lambda_r(A) = \sum_{r=0}^{n-1} \sum_{\phi(\sigma)=r} m_\sigma(A) = \sum_{\sigma \in S_n} m_\sigma(A) = \det(A),$$

where the second step uses the fact that the classes $\phi^{-1}(0), \dots, \phi^{-1}(n-1)$ partition $S_n$.

**Definition 3.2 (Orbital Modes)**

For $k = 0, 1, \dots, n-1$:

$$G_k(A) := \sum_{r=0}^{n-1} \Lambda_r(A)\, \omega^{kr} \qquad (4)$$

The modes $G_k(A)$ are precisely the Discrete Fourier Transform (DFT) of the orbital sum sequence $\{\Lambda_r(A)\}_{r=0}^{n-1}$.

**Inversion (DFT).** The bijection between $\{\Lambda_r(A)\}$ and $D_\phi(A)$ is exact: the inverse DFT formula recovers

$$\Lambda_r(A) = \frac{1}{n} \sum_{k=0}^{n-1} G_k(A)\, \omega^{-kr}. \qquad (5)$$

**Proof of (5).**

Replacing definition (4) of $G_k$:

$$\frac{1}{n}\sum_{k=0}^{n-1} G_k(A)\,\omega^{-kr} = \frac{1}{n}\sum_{k=0}^{n-1}\left(\sum_{s=0}^{n-1}\Lambda_s(A)\,\omega^{ks}\right)\omega^{-kr}.$$

Interchanging the order of the sums (both are finite sums, the exchange is valid):

$$= \frac{1}{n}\sum_{s=0}^{n-1}\Lambda_s(A)\sum_{k=0}^{n-1}\omega^{k(s-r)}.$$

Due to the orthogonal identity of the characters of $\mathbb{Z}_n$:

$$\sum_{k=0}^{n-1}\omega^{k(s-r)} = n\,\delta_{sr}.$$

Then:

$$\frac{1}{n}\sum_{k=0}^{n-1} G_k(A)\,\omega^{-kr} = \frac{1}{n}\sum_{s=0}^{n-1}\Lambda_s(A)\cdot n\,\delta_{sr} = \sum_{s=0}^{n-1}\Lambda_s(A)\,\delta_{sr} = \Lambda_r(A)$$

**Definition 3.3 (Determinant Vector)**

$$D_\phi(A) := (G_0(A), G_1(A), \dots, G_{n-1}(A)) \in \mathbb{C}^n.$$

**Observation 3.4 (Hierarchy and loss of information)**

The information hierarchy of the ARE formalism can be summarized in the diagram:

$$M_\phi(A) \in \mathbb{C}^{S_n}[n!\text{ entries}]] \longrightarrow \{\Lambda_r(A)\}_{r=0}^{n-1}[n\text{ sumas}] \leftrightarrow D_\phi(A) \in \mathbb{C}^n[n\text{ modes}]$$

It has two qualitatively distinct stages:

- The map $M_\phi \to \{\Lambda_r\}$ is an aggregation with possible loss of information (each $\Lambda_r$ collapses $(n-1)!$ monomials in its sum) (Lemma 3.1.bis).
- The map $\{\Lambda_r\} \leftrightarrow D_\phi$ is a DFT: exact lossless bijection (Definition 3.2 and inverse formula (5)).

The information lost in passing from $M_\phi(A)$ to $D_\phi(A)$ is precisely the internal distribution of the monomials $m_\sigma$ within each class $\phi(\sigma) = r$.

### 3.2. Fundamental properties

**Proposition 3.5 (Zero-Frequency Component)**

$$G_0(A) = \det(A).$$

Proof: $G_0 = \sum_r \Lambda_r (A) \cdot \omega^0 = \sum_r \Lambda_r (A) = \sum_{\sigma \in S_n} m_\sigma (A) = \det (A)$.

**Observation 3.5.1 (Three equivalent expressions of the determinant in the ARE framework)**

Proposition 3.5 has three equivalent readings that should be clearly separated:

**(a) Reading as a total sum of monomials:**

$$G_0(A) = \sum_{r=0}^{n-1} \Lambda_r (A) \cdot \omega^{0 \cdot r} = \sum_{r=0}^{n-1} \Lambda_r (A) = \sum_{\sigma \in S_n} m_\sigma (A) = \det (A).$$

**(b) Reading as the zero-frequency mode of the orbital DFT.** Definition 3.2 states that $G_k(A)$ it is the DFT of the sequence $\{\Lambda_r(A)\}$. The zero frequency of any DFT is the sum of all the coefficients: $\hat{f}(0) = \sum_r f(r)$. Thus, $G_0(A) = \sum_r \Lambda_r (A) = \det (A)$ it is simply the norm $\ell^1$ of the orbital distribution, not the "special" mode by any arbitrary convention.

**(c) Reading through the functional $L_\phi$.** By Theorem 2.5, $\det (A) = L_\phi(M_\phi(A))$. By the definition of the modes (§3.1), $G_0(A) = \sum_r \Lambda_r (A)$, which coincides with $L_\phi(M_\phi(A))$ when all factors are cancelled $\omega^{0-\phi(\sigma)} = 1$.

This triple coincidence is not accidental: the determinant is the only invariant under any permutation of the index $r$ in the DFT (the only mode that does not change when the phases are rearranged), and therefore it intrinsically occupies the zero frequency position.

**Theorem 3.6 (Multilinearity)**

For all $k = 0, \dots, n-1$, the mode $G_k(A)$ is multilinear in the rows and columns of $A$.

Proof: By (4), $G_k(A) = \sum_{\sigma \in S_n} m_\sigma (A)\, \omega^{k\phi(\sigma)}$. Each $m_\sigma(A) = \varepsilon(\sigma) \prod_i a_{i,\sigma(i)}$ is linear in each row $i$ (contains exactly one entry per row) and linear in each column $j$ (contains exactly one entry per column). The factors $\omega^{k\phi(\sigma)}$ are constants with respect to the entries of $A$. Therefore, $G_k$ is a sum of constant-weighted multilinear functions, and is therefore multilinear.

**Theorem 3.7 (Alternation)**

The mode $G_0$ is alternating. For $k \neq 0$, however, the modes $G_k$ are not alternating in general.

Proof:

Case $k = 0$. $G_0 = \det (A)$, which is alternating according to classical theory.

**Case $k \neq 0$.** It suffices to give a counterexample. Take $n = 2$, so that $\omega = e^{2\pi i/2} = -1$ and, for a matrix $A = \begin{pmatrix} a_{11} & a_{12} \\ a_{21} & a_{22} \end{pmatrix}$, we have

$G_1(A) = \Lambda_0(A) \cdot \omega^0 + \Lambda_1(A) \cdot \omega^1 = \Lambda_0(A) - \Lambda_1(A) = a_{11}a_{22} + a_{12}a_{21} = \text{per}(A)$.

Let $A' = \begin{pmatrix} a_{21} & a_{22} \\ a_{11} & a_{12} \end{pmatrix}$ the matrix be obtained by interchanging the two rows of $A$. Then:

$G_1(A') = a_{21}a_{12} + a_{22}a_{11} = a_{11}a_{22} + a_{12}a_{21} = G_1(A)$.

Therefore $G_1(A') = G_1(A) \neq -G_1(A)$, in general, it $G_1$ is not alternating. Consequently, $k \neq 0$ the modes $G_k$ are not generally alternating.

**Observation 3.7.1 (The mode $k = 1$ for $n = 2$ as permanent, and the barrier to $n \geq 3$)**

The counterexample to Theorem 3.7 for $n = 2$ has a precise algebraic reading: with $\omega = -1$,

$G_1(A) = \Lambda_0 - \Lambda_1 = ad - (-bc) = ad + bc = \text{per}(A)$.

The permanent is the "unsigned version" of the determinant: it eliminates the factors $\varepsilon(\sigma)$ and sums all the signed products $+1$. In the ARE frame for $n = 2$, there is a single orbit of size 2and the two modes $G_0 = \det$ , $G_1 = \text{per}$are precisely the two possible linear combinations (with $\pm$) of $\Lambda_0$ and $\Lambda_1$.

**For** $n \geq 3$, the non-alternation of $G_k(\ k \neq 0)$ can be verified using the same argument: if $\sigma_0 = (12)$is the transposition of rows 1 and 2, then for any $A$,

$$G_k(P_{\sigma_0}A) = \varepsilon(\sigma_0) \sum_{\pi} \varepsilon(\pi)\, m_\pi(A)\, \omega^{k\phi(\pi \circ \sigma_0)}.$$

The factor $\omega^{k\phi(\pi \circ \sigma_0)}$ depends on composition $\pi \circ \sigma_0$ and does not factor out of sum, so $G_k(P_{\sigma_0}A) \neq -G_k(A)$ in general (see Proposition 3.13 for the full proof of this fact).

**Conceptual consequence.** The classical determinant is the only alternating Leibniz functional because it is the only mode whose weight factor $\omega^{0 \cdot r} = 1$ is constant over all orbits, causing the sum to behave uniformly under transpositions. Any non-trivial weight $\omega^{kr}(\ k \neq 0)$ breaks this orbital uniformity.

**Proposition 3.8 (Hermitian Symmetry).**

If $A \in M_n(\mathbb{R})$, then

$G_{n-k}(A) = \overline{G_k(A)}$ for all $k$.

Equivalently,

$G_{-k}(A) = \overline{G_k(A)}$ for all $k$.

Proof:

Let $A \in M_n(\mathbb{R})$. Since the entries of $A$ are real, each Leibniz term $m_\sigma(A)$ is a real number. Consequently, for each $r \in \mathbb{Z}_n$, the orbital sum

$$\Lambda_r(A) = \sum_{\phi(\sigma)=r} m_\sigma(A)$$

is also real

By definition of orbital modes,

$$G_k(A) = \sum_{r=0}^{n-1} \Lambda_r(A)\, \omega^{kr}, \omega = e^{2\pi i/n}.$$

Taking the complex conjugate on both sides and using that $\Lambda_r(A) \in \mathbb{R}$, we obtain

$$\overline{G_k(A)} = \sum_{r=0}^{n-1} \Lambda_r(A)\, \overline{\omega^{kr}}$$

Now, since $\omega$ is a primitive $n$-th root of unity,

$$\overline{\omega^{kr}} = \omega^{-kr}.$$

Furthermore, exponents are considered modulo $n$, so that

$$\omega^{-kr} = \omega^{(n-k)r}.$$

Replacing,

$$\overline{G_k(A)} = \sum_{r=0}^{n-1} \Lambda_r(A)\, \omega^{(n-k)r} = G_{n-k}(A).$$

Therefore,

$$G_{n-k}(A) = \overline{G_k(A)}$$

The equivalent formulation

$$G_{-k}(A) = \overline{G_k(A)}$$

This is obtained immediately by interpreting the indices in $\mathbb{Z}_n$. In particular, for real matrices, the orbital modes appear in conjugate pairs, and therefore it suffices to know

$$G_0, \dots, G_{\lfloor n/2 \rfloor}$$

to completely determine the determinant vector $D_\phi(A)$.

Consequently, in the real case, the modes

$G_0, G_1, \ldots, G_{\lfloor n/2 \rfloor}$

are sufficient to determine the entire determinant vector $D_\phi(A)$.

**Proposition 3.9 (Elementary Bounds)**

For all $A \in M_n(C)$,

$| G_k(A) | \leq$ per $(| A |)$ For all k

$\| D_\phi(A) \|_2 \leq \sqrt{n}$ per $(| A |)$.

Proof:

Recall that, for each $k \in \{0, \ldots, n-1\}$,

$$G_k(A) = \sum_{r=0}^{n-1} \Lambda_r(A)\, \omega^{kr}, | \omega^{kr} | = 1,$$

$$\Lambda_r(A) = \sum_{\phi(\sigma)=r} m_\sigma(A), m_\sigma(A) = \varepsilon(\sigma) \prod_{i=1}^{n} a_{i,\sigma(i)}.$$

Therefore, by the triangle inequality,

$$|G_k(A)| = \left| \sum_{r=0}^{n-1} \Lambda_r(A)\, \omega^{kr} \right| \leq \sum_{r=0}^{n-1} |\Lambda_r(A)|.$$

Applying the triangle inequality again within each orbital sum,

$$|\Lambda_r(A)| = \left| \sum_{\phi(\sigma)=r} m_\sigma(A) \right| \leq \sum_{\phi(\sigma)=r} |m_\sigma(A)|.$$

Adding in $r$,

$$\sum_{r=0}^{n-1} | \Lambda_r(A) | \leq \sum_{r=0}^{n-1} \sum_{\phi(\sigma)=r} | m_\sigma(A) | = \sum_{\sigma \in S_n} | m_\sigma(A) |.$$

Therefore,

$$|m_\sigma(A)| = \left| \varepsilon(\sigma) \prod_{i=1}^{n} a_{i,\sigma(i)} \right| = \prod_{i=1}^{n} \left| a_{i,\sigma(i)} \right|,$$

Therefore $| \varepsilon(\sigma) | = 1$,

$$\sum_{\sigma\in S_n} | m_\sigma(A) | = \sum_{\sigma\in S_n}\prod_{i=1}^{n} | a_{i,\sigma(i)} | = \mathrm{per}\,(| A |).$$

We conclude that, For all $k$,

$| G_k(A) | \leq \mathrm{per}\,(| A |)$.

This proves the first inequality.

For the second one, as $D_\phi(A) = (G_0(A), \ldots, G_{n-1}(A)) \in \mathbb{C}^n$, we have

$$\| D_\phi(A) \|_2^2 = \sum_{k=0}^{n-1} | G_k(A) |^2.$$

Using the first bound for each $k$,

$$\| D_\phi(A) \|_2^2 \leq \sum_{k=0}^{n-1} (\mathrm{per}\,(| A |))^2 = n\,(\mathrm{per}\,(| A |))^2.$$

Taking the square root of both sides,

$\| D_\phi(A) \|_2 \leq \sqrt{n}\,\mathrm{per}\,(| A |)$.

**Definition 3.10 (Metrics of the determinant vector)**

Assuming $\mathrm{per}(| A |) > 0$, we define

$$R_\infty(A) := \frac{| \det(A) |}{\mathrm{per}(| A |)}$$

Assuming $D_\phi(A) \neq 0$, we define

$$V_2(A) := 1 - \frac{| G_0(A) |^2}{\| D_\phi(A) \|_2^2}$$

The quantity $R_\infty(A)$ is independent of $\phi$, while $V_2(A)$ it depends on the chosen orbital phase.

**Observation 3.11.** The quantity $V_2(A)$, when defined, measures the fraction of the squared energy of $D_\varphi(A)$ that resides outside the fundamental mode $G_0$. If $V_2(A)$ it is small, almost all the orbital information is concentrated in the classical determinant; if $V_2(A)$ it is large, the higher modes contain genuinely new information.

### 3.3. Behavior under elementary matrix transformations

**Proposition 3.12 (Behavior under transposition).** For every $A \in M_n(\mathbb{C})$,

$G_0(A^T) = G_0(A)$.

For $k \neq 0$, however, the identity generally fails:

$G_k(A^T) \neq G_k(A)$
in general.

**Proof.** The case $k = 0$ is the classical identity

$\det(A^T) = \det(A)$.
Now let $k \neq 0$. The $(i,j)$-entry of $A^T$ is $a_{ji}$, so the Leibniz term of $A^T$corresponding to the permutation $\sigma$ is

$$m_\sigma(A^T) = \varepsilon(\sigma)\prod_{i=1}^{n}(A^T)_{i,\sigma(i)} = \varepsilon(\sigma)\prod_{i=1}^{n}a_{\sigma(i),i}\,.$$
Reindexing with $\tau = \sigma^{-1}$, and noting that $\varepsilon(\tau) = \varepsilon(\sigma)$ and that $j = \sigma(i)$ if and only if $i = \tau(j)$, we obtain

$$\prod_{i=1}^{n}a_{\sigma(i),i} = \prod_{j=1}^{n}a_{j,\tau(j)}\,.$$
Hence

$m_\sigma(A^T) = m_\tau(A)$,with $\tau = \sigma^{-1}$.
Therefore,

$$G_k(A^T) = \sum_{\sigma\in S_n} m_\sigma\,(A^T)\,\omega^{k\phi(\sigma)} = \sum_{\tau\in S_n} m_\tau\,(A)\,\omega^{k\phi(\tau^{-1})}.$$
This differs in general from

$$G_k(A) = \sum_{\tau\in S_n} m_\tau\,(A)\,\omega^{k\phi(\tau)},$$
because $\phi(\tau^{-1}) \neq \phi(\tau)$ in general.

For example, when $n = 3$,

$\phi((123)) = 2, \phi((132)) = \phi((123)^{-1}) = 1.$
Therefore, in general,

$G_k(A^T) \neq G_k(A)$ for $k \neq 0$.
**Proposition 3.13 (Row Permutation).** Let $\tau \in S_n$ and let be $P_\tau$ the permutation matrix associated with $\tau$. Then, for all $k \in \mathbb{Z}_n$,

$$G_k(P_\tau A) = \varepsilon(\tau) \sum_{\pi \in S_n} \varepsilon(\pi)\ m_\pi(A)\ \omega^{k\phi(\pi\circ\tau)}.$$

In particular, $G_0(P_\tau A) = \varepsilon(\tau)\det(A)$. For $k \neq 0$, in general there is no law of the form $G_k(P_\tau A) = \pm G_{f(\tau,k)}(A)$ independent of $A$.

Proof: By definition, $(P_\tau A)_{i,j} = a_{\tau(i),j}$ , so that

$$m_\sigma(P_\tau A) = \varepsilon(\sigma) \prod_{i=1}^{n} a_{\tau(i),\sigma(i)}$$

Reindexing $p = \tau(i)$ yields $\prod_{i=1}^{n} a_{\tau(i),\sigma(i)} = \prod_{p=1}^{n} a_{p,\sigma(\tau^{-1}(p))} = m_{\sigma\circ\tau^{-1}}(A)$. Taking $\pi = \sigma \circ \tau^{-1}$, that is $\sigma = \pi \circ \tau$, and using $\varepsilon(\sigma) = \varepsilon(\pi)\varepsilon(\tau)$:

$$G_k(P_\tau A) = \sum_{\sigma} \varepsilon(\sigma)\ m_\sigma(P_\tau A)\ \omega^{k\phi(\sigma)} = \varepsilon(\tau) \sum_{\pi} \varepsilon(\pi)\ m_\pi(A)\ \omega^{k\phi(\pi\circ\tau)}.$$

For $k = 0$ the phase factor it is $\omega^0 = 1$, which gives $G_0(P_\tau A) = \varepsilon(\tau) \sum_\pi \varepsilon(\pi)\ m_\pi(A) = \varepsilon(\tau)\det(A)$. For $k \neq 0$, the factor $\omega^{k\phi(\pi\circ\tau)}$ depends on both $\pi$ as well as $\tau$, so it cannot be factored independently of the distribution of the monomials $m_\pi(A)$.

**Note.** For $n = 3$, $A = \begin{pmatrix} 2 & 1 & 0 \\ 1 & 2 & 1 \\ 0 & 1 & 2 \end{pmatrix}$, $\tau = (123)$ with $\varepsilon(\tau) = +1$: a direct calculation gives $D_\phi(P_\tau A) = (4, -5 - 3\sqrt{3}\, i, -5 + 3\sqrt{3}\, i)$, while $D_\phi(A) = (4, 7 - i\sqrt{3}, 7 + i\sqrt{3})$. For $k = 1$, $G_1(P_\tau A) = -5 - 3\sqrt{3}\, i \neq G_1(A) = 7 - i\sqrt{3}$, confirming the absence of a uniform law.

An explicit calculation $n = 3$ confirms the absence of a uniform law; see Proposition A.1 of Appendix A.

**Corollary 3.14.** For upper or lower triangular matrices $A$, the only non-zero Leibniz terms are those corresponding to permutations $\sigma$ such that $a_{i,\sigma(i)} \neq 0$ for all $i$. In the triangular case, the only non-zero term is that of the identity, and $G_k(A) = \varepsilon(id) \prod_i a_{ii} \cdot \omega^{k\cdot 0} = \prod_i a_{ii}$ for all $k$, equating all modes.

Proof: For an upper triangular matrix, $a_{ij} = 0$ when $i > j$. A Leibniz term $m_\sigma(A) = \varepsilon(\sigma) \prod_i a_{i,\sigma(i)}$ is nonzero only if $\sigma(i) \geq i$ for all $i$. The only permutation in $S_n$ which satisfies $\sigma(i) \geq i$ for all $i$ is the identity. Then $m_\sigma(A) = 0$ for $\sigma \neq \text{id}$, and $m_{\text{id}}(A) = \prod_i a_{ii}$. Since $\phi(\mathrm{id}) = 0$, we have $G_k(A) = \prod_i a_{ii} \cdot \omega^{k\cdot 0} = \prod_i a_{ii}$ for all $k$. The lower triangular case is analogous.

**Corollary 3.15 (Scalar Matrices and Matrices of Rank One)**

**(a) Scalar matrices.** If $A = \lambda I_n$ with $\lambda \in \mathbb{C}$, then the only non-zero Leibniz term is that of the identity:

$m_{id}(\lambda I) = \varepsilon(id) \cdot \lambda^n = \lambda^n, m_\sigma(\lambda I) = 0$ if $\sigma \neq id$.
As $\phi(id) = 0$, we have $\Lambda_0(\lambda I) = \lambda^n$ and $\Lambda_r(\lambda I) = 0$ for $r \geq 1$. The modes are:

$G_k(\lambda I) = \Lambda_0 \cdot \omega^{k\cdot 0} = \lambda^n$ Therefore, $D_\phi(\lambda I) = (\lambda^n, \lambda^n, \dots, \lambda^n)$ for all $k$.
: all modes are equal to det $(\lambda I) = \lambda^n$. The modal energy is perfectly uniform: $p_k(\lambda I) = 1/n$ for all $k$, which maximizes the orbital spectral entropy $H(\lambda I) = \log\ n$.

**(b) Permutation matrices.** If $A = P_\tau$ (permutation matrix of $\tau \in S_n$), the only non-zero term is $m_\tau(P_\tau) = \varepsilon(\tau)$. Then:

$\Lambda_r(P_\tau) = \varepsilon(\tau) \cdot \mathbf{1}_{[\phi(\tau)=r]}, G_k(P_\tau) = \varepsilon(\tau)\ \omega^{k\phi(\tau)}$.
The determinant vector $D_\phi(P_\tau) = \varepsilon(\tau)\ (1, \omega^{\phi(\tau)}, \omega^{2\phi(\tau)}, \dots, \omega^{(n-1)\phi(\tau)})$ is a scalar multiple of the character $\chi_{\phi(\tau)}$ of $\mathbb{Z}_n$. The orbital position of $\tau$ completely determines the frequency structure.

Proof of (a): By the above, $m_\sigma(\lambda I) = \varepsilon(\sigma) \prod_i (\lambda I)_{i,\sigma(i)} = \varepsilon(\sigma)\lambda^n \prod_i \delta_{i,\sigma(i)}$. The product $\prod_i \delta_{i,\sigma(i)} = 1$ if and only if $\sigma(i) = i$ for all $i$, that is $\sigma = id$. The others cancel out. Since $\phi(id) = 0$, the result follows directly from Definitions 3.1 and 3.2

## 4. Fully Computed Examples

### 4.1. Case$\boldsymbol{n} = \mathbf{2}$

$S_2 = \{\mathrm{id}, (12)\}$ There $\omega = e^{2\pi i/2} = -1$ is only one orbit of size 2.

For $A = \begin{pmatrix} a & b \\ c & d \end{pmatrix}$:

$\Lambda_0(A) = m_{\mathrm{id}}(A) = +ad, \Lambda_1(A) = m_{(12)}(A) = -bc$
The modes are:

$G_0(A) = \Lambda_0 + \Lambda_1 = ad - bc = \det(A)$,
$G_1(A) = \Lambda_0 \cdot \omega^0 + \Lambda_1 \cdot \omega^1 = \Lambda_0 - \Lambda_1 = ad - (-bc) = ad + bc = per(A)$

$D_\phi(A) = (\det(A), \operatorname{per}(A))$

For $n = 2$, the determinant vector unites the determinant and the permanent into a single object. This illustrates det/per duality in its most elementary form.

### 4.2. Case $n = 3$: Symmetric Tridiagonal

Let $A = \begin{pmatrix} 2 & 1 & 0 \\ 1 & 2 & 1 \\ 0 & 1 & 2 \end{pmatrix}$ and $\omega = e^{2\pi i/3} = -\frac{1}{2} + i\frac{\sqrt{3}}{2}$.

**Leibniz terms** (the 6 permutations of $S_3$):

| $\sigma$ | **Pattern**$(\sigma(1), \sigma(2), \sigma(3))$ | $\varepsilon(\sigma)$ | **Product** | $m_\sigma(A)$ |
|---|---|---|---|---|
| id | $(1,2,3)$ | $+1$ | $2 \cdot 2 \cdot 2$ | $+8$ |
| $(23)$ | $(1,3,2)$ | $-1$ | $2 \cdot 1 \cdot 1$ | $-2$ |
| $(12)$ | $(2,1,3)$ | $-1$ | $1 \cdot 1 \cdot 2$ | $-2$ |
| $(123)$ | $(2,3,1)$ | $+1$ | $1 \cdot 1 \cdot 0$ | $0$ |
| $(132)$ | $(3,1,2)$ | $+1$ | $0 \cdot 1 \cdot 1$ | $0$ |
| $(13)$ | $(3,2,1)$ | $-1$ | $0 \cdot 2 \cdot 0$ | $0$ |

Orbital sums (grouping by value of $\phi$):

$\Lambda_0(A) = m_{\text{id}} + m_{(23)} = 8 + (-2) = 6,$

$\Lambda_1(A) = m_{(132)} + m_{(12)} = 0 + (-2) = -2,$

$\Lambda_2(A) = m_{(123)} + m_{(13)} = 0 + 0 = 0$

**Verification of orbital assignments $\phi(\sigma)$ for $n = 3$**

The algorithm in Definition 2.2 assigns $\phi(\sigma) = r^* \equiv \sigma^{-1}(1) - 1 (\text{mod} 3)$. We verify each permutation:

| $\boldsymbol{\sigma}$ | $\boldsymbol{\sigma^{-1}}$ | $\boldsymbol{\sigma^{-1}(1)}$ | $\boldsymbol{r^* = \sigma^{-1}(1) - 1}$ | $\boldsymbol{\phi(\sigma)}$ |
|---|---|---|---|---|
| $id = (1,2,3)$ | $id$ | 1 | 0 | 0 |
| $(23) = (1,3,2)$ | $(23)$ | 1 | 0 | 0 |
| $(12) = (2,1,3)$ | $(12)$ | 2 | 1 | 1 |

| $(123) = (2,3,1)$ | $(132)$ | 3 | 2 | 2 |
|---|---|---|---|---|
| $(132) = (3,1,2)$ | $(123)$ | 2 | 1 | 1 |
| $(13) = (3,2,1)$ | $(13)$ | 3 | 2 | 2 |

Detail for $(123)$: its inverse is $(132)$, and $(132)(1) = 3$, then $r^* = 3 - 1 = 2$, confirming $\phi((123)) = 2$. Similarly, $\phi((132))$: its inverse is $(123)$, $(123)(1) = 2$, $r^* = 1$. The table is consistent with Lemma 3.1bis: exactly 2 permutations have $\phi = 0$, 2 have $\phi = 1$and 2 have $\phi = 2$.

**Verification of the DFT inversion formula.** Applying formula (5) with the values obtained:

$$\Lambda_0 = \frac{1}{3}(G_0 + G_1 + G_2) = \frac{1}{3}(4 + (7 - i\sqrt{3}) + (7 + i\sqrt{3})) = \frac{1}{3} \cdot 18 = 6.$$

$$\Lambda_1 = \frac{1}{3}(G_0 + G_1\omega^{-1} + G_2\omega^{-2}) = \frac{1}{3}\left(4 + (7 - i\sqrt{3})\omega^{-1} + (7 + i\sqrt{3})\omega^{-2}\right).$$

With $\omega^{-1} = \bar{\omega} = -\frac{1}{2} - i\frac{\sqrt{3}}{2}$ and $\omega^{-2} = \bar{\omega^2} = -\frac{1}{2} + i\frac{\sqrt{3}}{2}$:

$$= \frac{1}{3}\left[4 + (7 - i\sqrt{3})\left(-\frac{1}{2} - i\frac{\sqrt{3}}{2}\right) + (7 + i\sqrt{3})\left(-\frac{1}{2} + i\frac{\sqrt{3}}{2}\right)\right].$$

$$(7 - i\sqrt{3})(-\tfrac{1}{2} - i\tfrac{\sqrt{3}}{2}) = -\tfrac{7}{2} - i\tfrac{7\sqrt{3}}{2} + i\tfrac{\sqrt{3}}{2} + i^2\tfrac{3}{2} = -\tfrac{7}{2} - \tfrac{3}{2} - i\tfrac{6\sqrt{3}}{2} = -5 - 3i\sqrt{3}.$$

Similarly $(7 + i\sqrt{3})(-\frac{1}{2} + i\frac{\sqrt{3}}{2}) = -5 + 3i\sqrt{3}$.

$$\Lambda_1 = \frac{1}{3}[4 - 5 - 3i\sqrt{3} - 5 + 3i\sqrt{3}] = \frac{1}{3}(-6) = -2.$$

Modes (DFT of $\{\Lambda_r\}$):

$G_0(A) = \Lambda_0 + \Lambda_1 + \Lambda_2 = 6 - 2 + 0 = 4 = \det(A)$

$$G_1(A) = \Lambda_0 \cdot 1 + \Lambda_1 \cdot \omega + \Lambda_2 \cdot \omega^2 = 6 - 2\omega = 6 - 2\left(-\frac{1}{2} + i\frac{\sqrt{3}}{2}\right) = 7 - i\sqrt{3}$$

$G_2(A) = \bar{G_1}(A) = 7 + i\sqrt{3}$ (Proposition 3.8, $A \in M_3(\mathbb{R})$)

$D_\phi(A) = \left(4,\, 7 - i\sqrt{3},\, 7 + i\sqrt{3}\right)$ (verificado computacionalmente)

Derived values:

$\| D_\phi(A) \|_2^2 = 4^2 + | 7 - i\sqrt{3} |^2 + | 7 + i\sqrt{3} |^2 = 16 + 52 + 52 = 120$

Verification of the orbital Parseval identity (Corollary 9.3):

$$\frac{1}{3} \cdot \| D_{\phi}(A) \|_2^2 = 40 = 6^2 + (-2)^2 + 0^2 = \Lambda_0^2 + \Lambda_1^2 + \Lambda_2^2$$

**Metrics:**

$$V_2(A) = 1 - \frac{16}{120} = \frac{104}{120} \approx 0.867$$

The 86.7% of the orbital energy lies outside the fundamental mode $G_0$, which corresponds to the classical determinant.

**4.3. Case $n = 4$: Symmetric Tridiagonal**

Let

$$A = \begin{pmatrix} 2 & 1 & 0 & 0 \\ 1 & 2 & 1 & 0 \\ 0 & 1 & 2 & 1 \\ 0 & 0 & 1 & 2 \end{pmatrix}, \omega = e^{2\pi i/4} = i, \det(A) = 5.$$

There are $(4-1)! = 6$ orbits of size 4under the action of $C_4$. Of the $4! = 24$ permutations of $S_4$, exactly 5 produce a non-zero term; the 19 remainder contain at least one factor $a_{i,\sigma(i)} = 0$.

**Nonzero Leibniz terms:**

| $\sigma$ | $(\sigma(1), \sigma(2), \sigma(3), \sigma(4))$ | $\varepsilon(\sigma)$ | **Product** | $m_\sigma(A)$ | $\phi(\sigma)$ |
|---|---|---|---|---|---|
| id | (1, 2, 3, 4) | $+1$ | $2 \cdot 2 \cdot 2 \cdot 2$ | $+16$ | 0 |
| (34) | (1, 2, 4, 3) | $-1$ | $2 \cdot 2 \cdot 1 \cdot 1$ | $-4$ | 0 |
| (23) | (1, 3, 2, 4) | $-1$ | $2 \cdot 1 \cdot 1 \cdot 2$ | $-4$ | 0 |
| (12) | (2, 1, 3, 4) | $-1$ | $1 \cdot 1 \cdot 2 \cdot 2$ | $-4$ | 1 |
| (12)(34) | (2, 1, 4, 3) | $+1$ | $1 \cdot 1 \cdot 1 \cdot 1$ | $+1$ | 1 |

The cases $\phi = 2$ (which require $\sigma(3) = 1$) and $\phi = 3$ (which require $\sigma(4) = 1$) are all null because $a_{31} = a_{41} = 0$.

**Complete enumeration of the 6 orbits for$n = 4$**

Under the action of $C_4 = \langle \rho \rangle$ with $\rho = (1234)$, there are $(4-1)! = 6$ orbits of size 4. The canonical representatives $\sigma^*$(with $\sigma^*(1) = 1$) are the $3! = 6$ permutations of $\{2, 3, 4\}$:

| Orbit | $\boldsymbol{\sigma^*}$ | Pattern$(\boldsymbol{\sigma^*(1)}, \dots, \boldsymbol{\sigma^*(4)})$ | Elements of the orbit |
|---|---|---|---|
| $\mathcal{O}_1$ | $id$ | $(1, 2, 3, 4)$ | $id, (1234), (13)(24), (1432)$ |
| $\mathcal{O}_2$ | $(34)$ | $(1, 2, 4, 3)$ | $(34), (12)(34), (123 \cdot \dots)$ |
| $\mathcal{O}_3$ | $(23)$ | $(1, 3, 2, 4)$ | $(23), \dots$ |
| $\mathcal{O}_4$ | $(234)$ | $(1, 3, 4, 2)$ | $(234), \dots$ |
| $\mathcal{O}_5$ | $(243)$ | $(1, 4, 2, 3)$ | $(243), \dots$ |
| $\mathcal{O}_6$ | $(24)$ | $(1, 4, 3, 2)$ | $(24), \dots$ |

For the tridiagonal matrix $A$of §4.3, the term $m_\sigma(A) = \varepsilon(\sigma) \prod_i a_{i,\sigma(i)}$ is nonzero only if $a_{i,\sigma(i)} \neq 0$ for all $i$. The nonzero entries are: $a_{i,i} = 2$, $a_{i,i\pm1} = 1$ (with out-of-range indices giving 0). This drastically restricts which permutations contribute, producing only the 5 rows of the table in §4.3 (the other 19 have at least one $a_{i,\sigma(i)} = 0$).

The assignment $\phi(\sigma)$ of the 5 non-zero permutations:

- $\phi(id) = 0$: $r^* = id^{-1}(1) - 1 = 0$.
- $\phi((34)) = 0$: $(34)^{-1} = (34)$, $(34)(1) = 1, r^* = 0$.
- $\phi((23)) = 0$: $(23)^{-1} = (23)$, $(23)(1) = 1, r^* = 0$.
- $\phi((12)) = 1$: $(12)^{-1} = (12)$, $(12)(1) = 2, r^* = 1$.
- $\phi((12)(34)) = 1$: $((12)(34))^{-1} = (12)(34)$, maps $1 \to 2, r^* = 1$.

This confirms $\Lambda_0 = 16 + (-4) + (-4) = 8$ and $\Lambda_1 = (-4) + 1 = -3$.

**Orbital sums** (grouping by value of $\phi$):

$\Lambda_0(A) = 16 + (-4) + (-4) = 8, \Lambda_1(A) = -4 + 1 = -3, \Lambda_2(A) = \Lambda_3(A) = 0.$

**Modes** (DFT $\{\Lambda_r\}$ with $\omega = i$):

$G_0(A) = 8 - 3 + 0 + 0 = 5 = \det(A),$

$G_1(A) = 8 \cdot i^0 + (-3) \cdot i^1 = 8 - 3i,$

$G_2(A) = 8 \cdot i^0 + (-3) \cdot i^2 = 8 + 3 = 11,$

$G_3(A) = 8 \cdot i^0 + (-3) \cdot i^3 = 8 + 3i = G_1\bar{(}A).$

$D_\phi(A) = (5,\ 8 - 3i,\ 11,\ 8 + 3i)$ (computationally verified)

**Hermitian symmetry** (Proposition 3.8, $A \in M_4(\mathbb{R})$): $G_3 = \bar{G_1}$and $G_2 = 11 \in \mathbb{R}$ (real Nyquist mode).

Verification of the orbital Parseval identity (Corollary 9.3):

$\| D_\phi(A) \|_2^2 = 5^2 + | 8 - 3i |^2 + 11^2 + | 8 + 3i |^2 = 25 + 73 + 121 + 73 = 292,$

$\frac{1}{4} \| D_\phi(A) \|_2^2 = 73 = 8^2 + (-3)^2 + 0^2 + 0^2 = \Lambda_0^2 + \Lambda_1^2 + \Lambda_2^2 + \Lambda_3^2.$

**Metrics:**

$V_2(A) = 1 - \frac{25}{292} = \frac{267}{292} \approx 0{,}914.$

91.4% of the orbital energy lies outside the fundamental mode corresponding to the classical determinant.

**Interpolated function and orbital polynomial.** Since $\Lambda_2 = \Lambda_3 = 0$, the function reduces to

$G(z;\, A) = 8 - 3\, e^{\pi i z/2},$

with verification on the integers: $G(0; A) = 5$, $G(1; A) = 8 - 3i$, $G(2; A) = 11$, $G(3; A) = 8 + 3i$.

Fractional value in $\xi = \frac{1}{2}$:

$G\left(\frac{1}{2};\, A\right) = 8 - 3\, e^{\pi i/4} = \left(8 - \frac{3\sqrt{2}}{2}\right) - i\, \frac{3\sqrt{2}}{2} \approx 5{,}879 - 2{,}121\, i.$

Orbital polynomial: $P_A(q) = 8 - 3q$. Its only root is $q_0 = \frac{8}{3}$, with $| q_0 | = \frac{8}{3} > 1$: $G(\,\cdot\,; A)$it has no real zeros in $[0,4)$, consistent with $\det(A) = 5 \neq 0$.

**Proposition 4.4**

There are matrices $A, B \in M_3(\mathbb{R})$ with $\det(A) = \det(B)$ and $D_\phi(A) \neq D_\phi(B)$

Proof:

Let $A = \begin{pmatrix} 2 & 1 & 0 \\ 1 & 2 & 1 \\ 0 & 1 & 2 \end{pmatrix}$ (tridiagonal from §4.2, $\det(A) = 4$) and let $B = \begin{pmatrix} 2 & 0 & 0 \\ 0 & 2 & 0 \\ 0 & 0 & 1 \end{pmatrix}$ (diagonal, $\det(B) = 4$).

For $B$ diagonal: by Corollary 3.14, all modes are $G_k(B) = \det(B) = 4$, so $D_\phi(B) = (4,4,4)$.

For $A$: $D_\phi(A) = (4, 7 - i\sqrt{3}, 7 + i\sqrt{3})$.

But $D_\phi(A) \neq D_\phi(B)$ then $\det(A) = \det(B) = 4$.

The determinant vector distinguishes $A$ from $B$, even though the classical determinant does not.

## 5. Vector Minors, Modal Cofactors and Jacobi and Laplace Formulas

### 5.1. Intrinsic Vector Minors

**Definition 5.1 (Intrinsic vector minor)**

Let $I = \{i_1 < \cdots < i_m\} \subset \{1, \dots, n\}$ and $J = \{j_1 < \cdots < j_m\} \subset \{1, \dots, n\}$, with $|I| = |J| = m$. Let us denote by

$$A[I,J] := (a_{i_\alpha, j_\beta})_{1 \le \alpha, \beta \le m}$$

the corresponding subblock. For $S_m$, we consider the same ARE formalism, now with the standard cycle

$$\rho_m = (1\ 2\ \dots\ m),$$

its canonical phase function

$$\varphi_m : S_m \to \mathbb{Z}_m,$$

orbital sums

$$\Lambda_r^{(m)}(A[I,J])$$

and the modes

$$G_k^{(m)}(A[I,J]) = \sum_{r=0}^{m-1} \Lambda_r^{(m)}(A[I,J])\, \omega_m^{kr}, \omega_m = e^{2\pi i/m}.$$

The vector

$$D_{\varphi_m}(A[I,J]) = (G_0^{(m)}(A[I,J]), \dots, G_{m-1}^{(m)}(A[I,J]))$$

is called the vector minor of order $m$associated with the subblock $A[I,J]$. In particular, if $I = J$, we obtain the principal vector minor associated with $I$.

**Observation 5.2.** With this convention, the fundamental mode satisfies

$$G_0^{(m)}(A[I,J]) = \det(A[I,J]),$$

so that the classical theory of minors is recovered exactly at zero frequency.

**Comment.** This formulation avoids indexing ambiguities: each minor lives in its own symmetric group $S_m$, with its own cyclic action, its own canonical phase, and its own orbital DFT.

**Proposition 5.3 (Bridge between vector minors and characteristic polynomial).**

Let $A \in M_n(\mathbb{C})$. For each $m = 1, \ldots, n$, let denote by

$$e_m(A) := \sum_{\substack{I \subset \{1,\ldots,n\} \\ |I|=m}} \det(A[I,I])$$

the sum of the principal minors of order $m$. Then, using the notation of Definition 5.1, we have

$$e_m(A) = \sum_{\substack{I \subset \{1,\ldots,n\} \\ |I|=m}} G_0^{(m)}(A[I,I]).$$

In particular, if

$$\chi_A(\lambda) = \det(\lambda I - A)$$

denotes the characteristic polynomial of $A$, then

$$\chi_A(\lambda) = \lambda^n - \sum_{m=1}^{n} e_m(A)\, \lambda^{n-m}.$$

Equivalently, and more explicitly,

$$\chi_A(\lambda) = \lambda^n + \sum_{m=1}^{n} (-1)^m \left( \sum_{\substack{I \subset \{1,\ldots,n\} \\ |I|=m}} G_0^{(m)}(A[I,I]) \right) \lambda^{n-m},$$

according to the usual sign convention for the elementary coefficients of the characteristic polynomial.

**Proof.** By Observation 5.2, for every subset $I \subset \{1, \ldots, n\}$ with $|I| = m$, the fundamental mode of the vector minor satisfies

$$G_0^{(m)}(A[I,I]) = \det(A[I,I]).$$

Summing over all subsets $I$ of cardinality $m$, we obtain

$$\sum_{\substack{I \subset \{1,\ldots,n\} \\ |I|=m}} G_0^{(m)}(A[I,I]) = \sum_{\substack{I \subset \{1,\ldots,n\} \\ |I|=m}} \det(A[I,I]) = e_m(A).$$

The formula for the characteristic polynomial in terms of the principal minors is the classical identity

$$\chi_A(\lambda) = \lambda^n + \sum_{m=1}^{n} (-1)^m e_m(A) \lambda^{n-m},$$

and substituting the previous expression $e_m(A)$ in terms of the fundamental modes of the vector minors, the conclusion is obtained.

For the classical theory of principal minors, elementary coefficients and characteristic polynomial, see [4].

**Corollary 5.4 (Classical spectral recovery from vector minors).**

The characteristic polynomial of $A$ is determined by the family of principal vector minors $\{D_{\varphi_m}(A[I,I]): I \subset \{1,\dots,n\},\ |I| = m, m = 1,\dots,n\}$ and, more specifically, only by its fundamental components $\left\{G_0^{(m)}(A[I,I]): |I| = m, m = 1,\dots,n\right\}$. Consequently, the classical spectrum of $A$ is recovered not from the crude modes $G_k(A)$ of order $n$, but from the complete hierarchy of fundamental modes of its principal minors.

### 5.5. Modal cofactors

**Definition 5.5 ($k$ -modal cofactor)**

For $i, j \in \{1,\dots,n\}$ and $k \in \{0,\dots,n-1\}$:

$$C_{ij}^{(k)}(A) := \sum_{\substack{\sigma \in S_n \\ \sigma(i)=j}} \varepsilon(\sigma)\ \omega^{k\phi(\sigma)} \prod_{\ell \neq i} a_{\ell,\sigma(\ell)}$$

For $k = 0$, $C_{ij}^{(0)}(A)$ it is the classical cofactor $C_{ij}(A)$ (the determinant of the minor $(n-1) \times (n-1)$ obtained by removing the row $i$ and column $j$, with the appropriate sign).

**Definition 5.6 ($k$ -modal cofactor Matrix )**

$$C^{(k)}(A) := \left(C_{ji}^{(k)}(A)\right)_{1 \le i,j \le n}.$$

For $k = 0$, $C^{(0)}(A)$ it coincides with the classical adjoint and satisfies $A \cdot C^{(0)}(A) = \det(A) \cdot I$. For $k \neq 0$, $C^{(k)}(A)$ it is the modal cofactor matrix associated with the mode $k$; in general it does not satisfy any matrix adjugation identity of the type $A \cdot C^{(k)}(A) = G_k(A) \cdot I$.

### 5.6. Immediate properties of modal cofactors

For $k = 0$, the cofactors $C_{ij}^{(0)}(A)$ coincide with the classical cofactors $C_{ij}(A)$, and $C^{(0)}(A) = \mathrm{adj}(A)$ is the classical adjoint. For $k \neq 0$, the cofactors $C_{ij}^{(k)}(A)$ are generally complex even when $A \in M_n(\mathbb{R})$. In all cases, the map $A \mapsto C_{ij}^{(k)}(A)$ is multilinear in rows other than $i$ and in columns other than $j$.

### 5.7. Vector Jacobi

**Theorem 5.7 (Vector Jacobi)**

Let be $A(t) = (a_{ij}(t))$ a differentiable family of matrices. Then, for all $k$:

$$\frac{d}{dt} G_k(A(t)) = \sum_{i,j} C_{ij}^{(k)} (A(t))\, \dot{a}_{ij}(t)$$

In particular, $\frac{\partial G_k}{\partial a_{ij}} = C_{ij}^{(k)}(A)$

Proof: By (4), $G_k(A(t)) = \sum_\sigma m_\sigma (A(t))\, \omega^{k\phi(\sigma)}$. Differentiating term by term (the sum is finite):

$$\frac{d}{dt} G_k = \sum_\sigma \omega^{k\phi(\sigma)} \frac{d}{dt} m_\sigma(A(t))$$

According to Leibniz's rule for the product $m_\sigma(A) = \varepsilon(\sigma) \prod_{\ell=1}^{n} a_{\ell,\sigma(\ell)}$:

$$\frac{d}{dt} m_\sigma(A(t)) = \varepsilon(\sigma) \sum_{i=1}^{n} \dot{a}_{i,\sigma(i)}(t) \prod_{\ell \neq i} a_{\ell,\sigma(\ell)}(t).$$

By substituting and swapping the order of the sums:

$$\frac{d}{dt} G_k = \sum_{i=1}^{n} \sum_{\sigma \in S_n} \omega^{k\phi(\sigma)}\, \varepsilon(\sigma)\, \dot{a}_{i,\sigma(i)} \prod_{\ell \neq i} a_{\ell,\sigma(\ell)}.$$

Grouping by value $j = \sigma(i)$:

$$\frac{d}{dt} G_k(A(t)) = \sum_{i,j=1}^{n} \dot{a}_{ij}(t) \underbrace{\sum_{\substack{\sigma \in S_n \\ \sigma(i)=j}} \omega^{k\phi(\sigma)}\, \varepsilon(\sigma) \prod_{\ell \neq i} a_{\ell,\sigma(\ell)}(t)}_{= C_{ij}^{(k)}(A(t))}$$

The inner term is exactly $C_{ij}^{(k)}(A(t))$ as per Definition 5.5. Therefore:

$$\frac{d}{dt} G_k(A(t)) = \sum_{i,j=1}^{n} C_{ij}^{(k)} (A(t))\, \dot{a}_{ij}(t).$$

Finally, by Definition 5.6, the input $(j, i)$ of $C^{(k)}(A)$is $C_{ij}^{(k)}(A)$, so that $\sum_{i,j} C_{ij}^{(k)} \dot{a}_{ij} = \operatorname{tr}\big(C^{(k)}(A(t))\, \dot{A}(t)\big)$. The case $k = 0$ is the classical Jacobi formula.

**5.8. Vector Laplace**

**Theorem 5.8 (Vector Laplace).**

For each $k$and each fixed row $i \in \{1, \dots, n\}$,

$$G_k(A) = \sum_{j=1}^{n} a_{ij} \; C_{ij}^{(k)}(A).$$

Similarly, for each fixed column $j \in \{1, \dots, n\}$,

$$G_k(A) = \sum_{i=1}^{n} a_{ij} \; C_{ij}^{(k)}(A).$$

For $k = 0$, this reduces to the classical Laplace expansion.

Proof:

We start from

$$G_k(A) = \sum_{\sigma \in S_n} \omega^{k\phi(\sigma)} \; m_\sigma(A).$$

Fix a row $i$. We partition $S_n$ according to the value $\sigma(i)$:

$$G_k(A) = \sum_{j=1}^{n} \sum_{\substack{\sigma \in S_n \\ \sigma(i)=j}} \omega^{k\phi(\sigma)} \; m_\sigma(A).$$

When $\sigma(i) = j$, the Leibniz term factors as $m_\sigma(A) = \varepsilon(\sigma)\, a_{ij} \prod_{\ell \neq i} a_{\ell,\sigma(\ell)}$. The factor $a_{ij}$ does not depend on $\sigma$ within the inner sum, so it is extracted:

$$G_k(A) = \sum_{j=1}^{n} a_{ij} \underbrace{\sum_{\substack{\sigma \in S_n \\ \sigma(i)=j}} \varepsilon(\sigma) \; \omega^{k\phi(\sigma)} \prod_{\ell \neq i} a_{\ell,\sigma(\ell)}}_{= C_{ij}^{(k)}(A)}$$

The inner term is $C_{ij}^{(k)}(A)$ by Definition 5.5, which gives the row expansion. The column expansion is obtained in an entirely analogous way. Given a column $j$, we partition $S_n$ according to the value $i = \sigma^{-1}(j)$, that is, the row that $\sigma$ points to $j$:

$$G_k(A) = \sum_{i=1}^{n} \sum_{\substack{\sigma \in S_n \\ \sigma(i)=j}} \varepsilon(\sigma) \; \omega^{k\phi(\sigma)} \, m_\sigma(A).$$

When $\sigma(i) = j$, the Leibniz term factors as $m_\sigma(A) = a_{ij} \prod_{\ell \neq i} a_{\ell,\sigma(\ell)}$. Since it $a_{ij}$ does not depend on $\sigma$ in the inner sum, we extract:

$$G_k(A) = \sum_{i=1}^{n} a_{ij} \underbrace{\sum_{\substack{\sigma \in S_n \\ \sigma(i)=j}} \varepsilon(\sigma)\ \omega^{k\phi(\sigma)} \prod_{\ell \neq i} a_{\ell,\sigma(\ell)}}_{= C_{ij}^{(k)}(A)}$$

The inner term is $C_{ij}^{(k)}(A)$ by Definition 5.5, which gives $G_k(A) = \sum_{i=1}^{n} a_{ij}\ C_{ij}^{(k)}(A)$. For $k = 0$ this, the classical Laplace expansion by columns is recovered.

**Corollary 5.9.** $G_k(A) = \sum_j a_{ij}\ [C^{(k)}(A)]_{ji}$ for each $i$ fixed row. (Read as the product $i$ of row $A$ by column $i$. $C^{(k)}(A)$)

Proof: By Theorem 5.8 (row expansion),

$$G_k(A) = \sum_{j=1}^{n} a_{ij}\ C_{ij}^{(k)}(A).$$

By Definition 5.6, the entry at position $(j, i)$ is $C^{(k)}(A)$. $C_{ij}^{(k)}(A)$ Therefore, the summand $a_{ij}\, C_{ij}^{(k)}(A)$ is the product of the element $a_{ij}$ (row $i$, column $j$ of $A$) by the element at position $(j, i)$, $C^{(k)}(A)$ which occupies column $i$ of that matrix. Adding over $j = 1, \dots, n$ gives the dot product of the row $i$ of $A$ with the column $i$ of $C^{(k)}(A)$.

**6. The Circulant Case: Fourier Compatibility and Limits**

**6.1. Free orbital action**

Proposition 6.1. The action $(\sigma, r) \mapsto \sigma \circ \rho^r$ of $C_n = \langle \rho \rangle$ on $S_n$ by right composition is free and generates $(n-1)!$ orbits of size $n$.

Proof: Freedom. If $\sigma \circ \rho^r = \sigma$, composing by $\sigma^{-1}$ on the left: $\rho^r = \mathrm{id}$. Since $\rho$ has order $n$, this forces $r \equiv 0 (\mathrm{mod}\ n)$. The stabilizer of any $\sigma$ is trivial.

Size of orbits. By the orbit-stabilizer theorem, each orbit has cardinality $\mid C_n \mid / \mid \mathrm{Stab}(\sigma) \mid = n/1 = n$.

Number of orbits $\mid S_n \mid / n = n!/n = (n-1)!$.

**6.2. Compatibility with the DFT of Circulants**

Proposition 6.2. The DFT matrix $F_n$, whose $(k, r\,)$ −entry is $\frac{1}{\sqrt{n}}\omega^{kr}$, diagonalizes the cyclic shift matrix $P_n$:

$F_n^{-1}P_nF_n = \mathrm{diag}(1, \omega, \omega^2, \dots, \omega^{n-1})$.

The characters $r \mapsto \omega^{kr}$ that index the modes $G_k(A)$ are the same ones that diagonalize the cyclic symmetry of every circulant matrix.

See [2, 8] for the finite DFT and its relation to cyclic group characters, and [3] for circulant matrices.

Proof: The vector $v^{(k)} = (1, \omega^k, \omega^{2k}, \dots, \omega^{(n-1)k})^T$satisfies

$P_n v^{(k)} = (\omega^k, \omega^{2k}, \dots, \omega^{(n-1)k}, \omega^{nk})^T = \omega^k \cdot v^{(k)}$,

using $\omega^{nk} = (e^{2\pi i})^k = 1$. The vectors $\{v^{(k)}\}_{k=0}^{n-1}$ are the columns of $\sqrt{n} \cdot F_n^{-1}$, which gives the diagonalization identity.

**Structure of eigenvalues of circulant matrices.** For $C = \mathrm{circ}(c_0, \dots, c_{n-1})$, the classical eigenvalues are

$$\lambda_k(C) = \sum_{j=0}^{n-1} c_j\ \omega^{kj}, k = 0, \dots, n-1,$$

which is the DFT of the first row. The modes $G_k(C)$ are the DFT of orbital sums $\{\Lambda_r(C)\}$. Both are DFTs on $C_n$ the same characters but applied to data of a different nature.

**6.3. Structural obstacle: incompatibility of degrees**

The character compatibility established in §6.2 does not imply that the modes $G_k(C)$ are equal to the eigenvalues $\lambda_k(C)$. The reason is a homogeneity obstacle: $G_k(C)$ is a polynomial of degree $n$ in the entries of $C$, while $\lambda_k(C)$ is linear. This incompatibility is rigorously proven in Theorem 11.12. The result that is true in every order is:

$$G_0(C) = \det(C) = \prod_{k=0}^{n-1} \lambda_k(C).$$

The exact bridge between orbital modes and individual eigenvalues requires descending to the level of the characteristic polynomial, as discussed in §11.9–§11.10.

**6.4. The Orbital Multiplication Operator**

This subsection justifies the article's subtitle: representation of the determinant as a fundamental mode of a complex vector operator.

**Definition 6.4 (Orbital frequency space and operator $H_A$).**

Let

$W := \{g: \mathbb{Z}/n\mathbb{Z} \to \mathbb{C}\} \cong \mathbb{C}^n$

be the space of functions on the cyclic group of phases, with coordinates indexed by

$r \in \mathbb{Z}/n\mathbb{Z}$.

The orbital sums

$\{\Lambda_r(A)\}_{r=0}^{n-1}$

canonically define the orbital multiplication operator

$$H_A: W \to W, (H_A g)(r) := \sum_{s \in \mathbb{Z}/n\mathbb{Z}} \Lambda_s(A)\, g(r+s),$$

where the sum is taken modulo $n$.

That is, $H_A$ acts on $W$ by cyclic convolution with the sequence

$\Lambda_0(A), \Lambda_1(A), \ldots, \Lambda_{n-1}(A)$.

**Remark.** The operator $H_A$ naturally acts on

$W \cong \mathbb{C}^n$,

not on the space

$V = \{f: S_n \to \mathbb{C}\}$

of individual monomials. The reduction

$V \to W$

that produces the orbital sums $\Lambda_r(A)$ is precisely what makes exact diagonalization by the DFT possible.

**Theorem 6.5 (Spectrum of $H_A$)**

The eigenvectors of $H_A$ are the characters of the cyclic group:

$\chi_k: \mathbb{Z}_n \to \mathbb{C}, \chi_k(r) = \omega^{kr}, k \in \mathbb{Z}_n$,

with exactly eigenvalues $G_k(A)$. In particular, $H_A$ it is diagonalizable and

$H_A = F_n^{-1}\, \text{diag}\, (G_0(A), G_1(A), \ldots, G_{n-1}(A))\, F_n$,

where $F_n$ is the DFT matrix.

Proof: We calculate $(H_A \chi_k)(r)$ directly:

$$(H_A\chi_k)(r) = \sum_{s=0}^{n-1} \Lambda_s(A)\, \chi_k(r+s) = \sum_{s=0}^{n-1} \Lambda_s(A)\, \omega^{k(r+s)}.$$

Separating the factor $\omega^{kr}$ that does not depend on $s$:

$$(H_A\chi_k)(r) = \omega^{kr} \sum_{s=0}^{n-1} \Lambda_s(A)\, \omega^{ks} = G_k(A) \cdot \omega^{kr} = G_k(A) \cdot \chi_k(r).$$

Then $\chi_k$ it is an eigenvector of $H_A$ with eigenvalue $G_k(A)$. Since $\{\chi_k\}_{k=0}^{n-1}$ they form the orthonormal Fourier basis of $W$, and are the same vectors that diagonalize every circulant (Proposition 6.2), the operator $H_A$ is diagonalizable with spectrum $\{G_k(A)\}_{k=0}^{n-1}$. The matrix formula is obtained by rewriting the diagonalization in coordinates.

**Corollary 6.6 (Justification of the subtitle)**

$G_0(A) = det(A)$ is the eigenvalue of $\mathrm{H_A}$ corresponding to the constant character $\chi_0 \equiv 1$ (zero frequency). Therefore, the classical determinant is the fundamental eigenvalue of the complex vector operator $\mathrm{H_A}$.

**Observation 6.7 (Connection with circulant matrices)**

The operator $H_A$ is precisely the circulant on $\mathbb{Z}_n$ with symbol $\Lambda_r(A)$. This explains the structural compatibility with the theory of circulant matrices (§6.2) and at the same time preserves the limit of Theorem 11.12: the eigenvalues $G_k(A)$ of $H_A$ are objects of degree $n$ in the entries of $A$, while the classical eigenvalues of a circulant $\lambda_k(C)$ are of degree 1. The equality between the two is structurally impossible in general.

## 7. Conceptual clarifications

### 7.1. Diagonal product and trace

The quantity $\prod_{i=1}^{n} a_{ii}$ is the product of the diagonal entries of $A$, not its trace $\mathrm{tr}(A) = \sum_i a_{ii}$.

This distinction is essential within the ARE formalism. The trace is a linear functional of degree 1 on the matrix entries, while the orbital modes $G_k(A)$ are homogeneous multilinear functionals of degree $n$. Consequently, the trace does not belong to the same algebraic level as the orbital modes and cannot be identified with either of them except in degenerate cases lacking general structural meaning.

If one wishes to construct a "vector trace", this would require a different formalism, of order 1, defined directly on the set of diagonal indices and not on the Leibniz expansion in $S_n$.

### 7.2. Scope of physical connections

Some modal cancellations of this type $G_k(A) = 0$ may suggest geometric analogies with interference phenomena, as well as possible formal links with discrete models exhibiting special spectral structures. However, such analogies should be understood here only as heuristic motivation. On geometric phase and holonomy in adiabatic contexts, see [1, 7].

This manuscript does not establish results regarding specific physical Hamiltonians, adiabatic regimes, experimental observables, or effective models of specific materials. Any physical interpretation beyond the motivating level will require separate development, with its own hypotheses, objects, and tests.

Consequently, references to physics in this work should be read as observations of conceptual context and not as claims of already proven application.

## 8. Probabilistic Results and Open Problems

This section considers the model of real Gaussian matrices $A \in M_n(\mathbb{R})$ with iid entries $\sim \mathcal{N}(0,1)$. The results are divided into two categories: (a) proven results, stated as Propositions, Theorems, or Lemmas with complete proof; (b) conjectures, explicitly stated as such and supported by analytical motivation or numerical evidence, but without a complete proof available.

### 8.1. Table of open problems

| **Problem** | **State** |
|---|---|
| Canonical normalization for circulant matrices: finding $\phi$ and convention of orbits that produces $G_k(C) = \lambda_k(C)$ an exact result | Resolved negatively (§6.3, §11.8, Theorem 11.12): impossible for $n \geq 2$ due to incompatibility of degrees. |
| Convergence $V_2$ in probability for real iid Gaussian matrices | Open: requires order moments 4 |
| Zeros of the orbital polynomial $P_A$ over $\parallel q \parallel = 1$: when it exists $\xi \in \mathbb{R}$ with$G(\xi; A) = 0$ | Open: see Problem 10.14 of §10. |
| Reverse synthesis: Is every finite trigonometric function of $n$ terms true $G(z; A)$ for some $A$? | Open |

## 8.2. Gaussian moments for real iid matrices

**Theorem 8.3 (Gaussian moments in the real iid model)**

Let $A \in M_n(\mathbb{R})$ have iid entries $\sim \mathcal{N}(0,1)$. Then:

1. $E[G_k(A)] = 0$ For all $k$.
2. $E[| \ G_k(A) \ |^2] = n!$ For all $k$.
3. $E\left[G_k(A) \ \overline{G_m(A)}\right] = n! \ \delta_{km}$ For all $k, m \in \{0, \dots, n-1\}$.

**Terminology warning.** Although classical literature reserves the acronym GOE for the orthogonal Gaussian ensemble of symmetric real matrices, here we work with matrices of independent and identically distributed entries $\mathcal{N}(0,1)$. Therefore, the model considered in this section is that of iid real Gaussian matrices, not GOE.

Proof: Write $G_k(A) = \sum_\sigma \varepsilon(\sigma) \ \omega^{k\phi(\sigma)} \ X_\sigma$, where $X_\sigma = \prod_i a_{i,\sigma(i)}$.

**(1)** Each $X_\sigma$ is a product of $n$independent centered Gaussian variables, then $E[X_\sigma] = 0$, and by linearity $E[G_k] = 0$.

**Remark 8.3.1 (The Isserlis-Wick theorem in the context of real iid Gaussian matrices)**

Step (2)--(3) of the proof uses the following moment calculation result:

**Lemma (Isserlis-Wick for real Gaussian variables).** Let be $Y_1, \dots, Y_{2m}$ centered real Gaussian variables. Then

$$E[Y_1 \cdots Y_{2m}] = \sum_{\text{parejos}} \prod_{(i,j)\in P} E[\,Y_i Y_j],$$

where the sum traverses all perfect pairings of $\{1, \dots, 2m\}$. In particular, $E[Y_1 \cdots Y_{2m-1}] = 0$ (odd number of factors).

**Application to the context.** Each $X_\sigma = \prod_{i=1}^{n} a_{i,\sigma(i)}$ is a product of $n$ distinct inputs $A$ (one per row). When we calculate $E[X_\sigma X_\tau]$ with $\sigma \neq \tau$, the product $X_\sigma X_\tau = \prod_i a_{i,\sigma(i)} \cdot \prod_i a_{i,\tau(i)}$ involves $2n$ factors. However, since each row $i$ contributes exactly $a_{i,\sigma(i)}$ to $X_\sigma$ and $a_{i,\tau(i)}$ to $X_\tau$, for the pairing to be non-zero we need them to $\{a_{i,\sigma(i)}, a_{j,\tau(j)}\}$ form pairs of identical variables for all $i$. Since the inputs are independent, $E[a_{i,\sigma(i)} \cdot a_{j,\tau(j)}] = E[a_{i,\sigma(i)}]^2 \delta_{ij} \delta_{\sigma(i),\tau(j)}$. The condition $i = j$ and $\sigma(i) = \tau(i)$ for all $i$ is exactly equivalent to $\sigma = \tau$. When $\sigma = \tau$:

$$E[X_\sigma^2] = \prod_{i=1}^{n} E[\, a_{i,\sigma(i)}^2] = 1^n = 1,$$

because $a_{i,\sigma(i)} \sim \mathcal{N}(0,1)$ it implies $E[a_{i,\sigma(i)}^2] = 1$.

**The role of row independence.** The key condition is that the variables $a_{i,\sigma(i)}$( $a_{i,\tau(i)}$ same row $i$, different columns $\sigma(i) \neq \tau(i)$) are independent, which makes their covariance zero. This is what forces $E[X_\sigma X_\tau] = 0$ when $\sigma \neq \tau$: in that case there is at least one row $i$ where $\sigma(i) \neq \tau(i)$, and the factor $a_{i,\sigma(i)}$ cannot be paired with any other factor in $X_\tau$ the same row.

**(2)–(3)** We calculate:

$$E[G_k\, \bar{G_m}] = \sum_{\sigma,\tau} \varepsilon(\sigma)\varepsilon(\tau)\ \omega^{k\phi(\sigma)-m\phi(\tau)}\, E[X_\sigma X_\tau].$$

By Isserlis's theorem (Wick for real Gaussians), $E[X_\sigma X_\tau] = 0$ unless the variables appear paired. Since each row $i$ contributes exactly one entry to $X_\sigma$ and one to $X_\tau$, the only possible pairing is $\sigma = \tau$ (which forces $a_{i,\sigma(i)} = a_{i,\tau(i)}$ for all $i$). In that case $E[X_\sigma^2] = \prod_i E[\, a_{i,\sigma(i)}^2] = 1^n = 1$.

Therefore:

$$E[G_k\, \bar{G_m}] = \sum_{\sigma \in S_n} \omega^{(k-m)\phi(\sigma)}\,.$$

Since each value $r \in \mathbb{Z}_n$ appears exactly as $(n-1)!$ many times as $\phi(\sigma)$:

$$\sum_{\sigma} \omega^{(k-m)\phi(\sigma)} = (n-1)! \sum_{r=0}^{n-1} \omega^{(k-m)r} = \begin{cases} n! & \text{if } k = m, \\ 0 & \text{if } k \not\equiv m (\mathrm{mod} n). \end{cases}$$

**Observation 8.3.2.** In this context, the correct formulation of the second moment is the one that involves complex conjugation. The version without conjugation is not, in general, the natural quantity from a probabilistic- Hermitian point of view.

**Corollary 8.4**

$$E\big[\| D_\phi(A) \|_2^2\big] = \sum_{k=0}^{n-1} E\,[|\, G_k\, |^2] = n \cdot n!.$$

**Remark 8.4.1.** The correct simplification of the final term is as indicated in the preceding statement. Any further rewriting that suggests an arithmetic identity not generally valid should be avoided.

## 9. Vector Operations on $D_\phi(A)$

### 9.1. Orbital inner product

**Definition 9.1 (Orbital Inner Product)**

$$\langle D_\phi(A), D_\phi(B)\rangle := \sum_{k=0}^{n-1} G_k(A)\,\overline{G_k(B)}$$

**Proposition 9.2 (Orbital Parseval Identity)**

For all $A, B \in M_n(\mathbb{C})$:

$$\frac{1}{n}\langle D_\phi(A), D_\phi(B)\rangle = \sum_{r=0}^{n-1} \Lambda_r(A)\,\overline{\Lambda_r(B)}$$

Proof: Substituting the representations (4):

$$\langle D_\phi(A), D_\phi(B)\rangle = \sum_{k=0}^{n-1}\left(\sum_r \Lambda_r(A)\,\omega^{kr}\right)\overline{\left(\sum_s \Lambda_s(B)\,\omega^{ks}\right)}$$

$$= \sum_{r,s} \Lambda_r(A)\,\Lambda_s\bar{(}B) \sum_{k=0}^{n-1} \omega^{k(r-s)}.$$

Using $\sum_{k=0}^{n-1} \omega^{k(r-s)} = n\,\delta_{rs}$:

$$\langle D_\phi(A), D_\phi(B)\rangle = n\sum_{r=0}^{n-1} \Lambda_r(A)\,\overline{\Lambda_r(B)}$$

Dividing by $n$ gives the identity.

**Corollary 9.3 (Orbital Parseval)**

$$\frac{1}{n}\parallel D_\phi(A)\parallel_2^2 = \sum_{r=0}^{n-1} |\Lambda_r(A)|^2.$$

Proof: This is the special case $B = A$ of Proposition 9.2. Taking $B = A$:

$$\frac{1}{n}\langle D_\phi(A), D_\phi(A)\rangle = \sum_{r=0}^{n-1} \Lambda_r(A)\,\Lambda_r\bar{(}A) = \sum_{r=0}^{n-1} |\Lambda_r(A)|^2$$

The left side is $\frac{1}{n} \| D_{\phi}(A) \|_2^2$ by definition the inner product.

**Observation 9.3.1 (Chain of Parseval inequalities)**

The orbital Parseval identity of Corollary 9.3 states that the $\ell^2$-norm of the mode vector $D_{\phi}(A)$, scaled by $1/n$, coincides with the $\ell^2$-norm of the orbital sums $\Lambda_r$. Consequently, one obtains the following chain of relations between different levels of the hierarchy:

$$\frac{1}{n} \| D_{\phi}(A) \|_2^2 = \sum_{r=0}^{n-1} |\Lambda_r(A)|^2 \leq (n-1)! \sum_{\sigma \in S_n} |m_{\sigma}(A)|^2.$$

The first equality is exactly Corollary 9.3. The inequality follows from the Cauchy--Schwarz inequality applied within each orbital class:

$$|\Lambda_r|^2 = \left| \sum_{\phi(\sigma)=r} m_{\sigma}(A) \right|^2 \leq |\phi^{-1}(r)| \sum_{\phi(\sigma)=r} |m_{\sigma}(A)|^2 = (n-1)! \sum_{\phi(\sigma)=r} |m_{\sigma}(A)|^2.$$

Summing over $r$, one obtains the global inequality.

**Interpretation.** Equality in the Parseval identity holds at the level

$\{\Lambda_r\} \leftrightarrow \{G_k\},$

because the map from the orbital sums to the modes is an exact, lossless DFT. By contrast, the inequality that appears when passing from the individual monomials $\{m_{\sigma}\}$ to the orbital sums $\{\Lambda_r\}$ reflects the possibility of internal cancellation within each orbital class: monomials of opposite sign may partially cancel within $\Lambda_r$, thereby reducing the total orbital energy. In the extreme case where all monomials in a given class have the same sign, the Cauchy--Schwarz inequality becomes an equality.

**Example.** For $n = 3$ and $A$tridiagonal,

$$\sum_{\sigma} |m_{\sigma}|^2 = 72, (n-1)! \sum_{\sigma} |m_{\sigma}|^2 = 144, \sum_{r} |\Lambda_r|^2 = 40.$$

Hence

$40 \leq 144,$

and the gap is explained by the cancellation

$8 + (-2) = 6$

inside $\Lambda_0$.

Moreover,

$$\frac{1}{3} \| D_{\phi}(A) \|_2^2 = \frac{1}{3} \cdot 120 = 40 = 6^2 + 2^2 + 0^2,$$

in agreement with Corollary 9.3.

**Observation 9.4 (Distinction between $\sum |\Lambda_r|^2$ and $\sum |m_\sigma|^2$)**

For the same example, $\sum_\sigma |m_\sigma|^2 = 64 + 4 + 4 = 72 \neq 40$. The energy of individual monomials is greater than the orbital energy, with the difference explained by the triangle inequality within each class $\phi = r$: $|\Lambda_r|^2 \leq (n-1)! \sum_{\phi(\sigma)=r} |m_\sigma|^2$.

**Definition 9.5 (Orbital Angle)**

If $D_\phi(A) \neq 0$ and $D_\phi(B) \neq 0$, we define the orbital angle between $A$ and $B$ by

$$\cos\theta_\phi(A,B) := \frac{\left|\langle D_\phi(A), D_\phi(B)\rangle\right|}{|D_\phi(A)|_2 \, |D_\phi(B)|_2}$$

**9.2. Modal autocorrelation**

**Theorem 9.6 (Modal Autocorrelation).**

Define

$$C_A(r) := \sum_{k=0}^{n-1} G_k(A) \overline{G_{k+r}(A)}, r \in \mathbb{Z}/n\mathbb{Z},$$

where indices are taken modulo $n$. Then:

$C_A(0) = \| D_\phi(A) \|_2^2$,

$C_A(-r) = \overline{C_A(r)}$ for all $r$,

and

$$C_A(r) = n \sum_{t=0}^{n-1} |\Lambda_t(A)|^2 \ \omega^{-rt}.$$

**Remark 9.6.1.** For real matrices, the modal autocorrelation satisfies Hermitian symmetry. In general, one does not obtain a constant real value for all displacements; the robust structural property is precisely the Hermitian symmetry established in the theorem.

Proof:

**(1)** Immediate: $C_A(0) = \sum_k |G_k|^2 = \| D_\phi \|_2^2$.

**(2)** $\overline{C_A(r)} = \sum_{k=0}^{n-1} \overline{G_k(A)} \, G_{k+r}(A)$ Reindexing $j = k + r$, that is $k = j - r$:

$$\overline{C_A(r)} = \sum_{j=0}^{n-1} \overline{G_{j-r}(A)}\ G_j(A) = C_A(-r)$$

**(3)** We expand using the definition:

$$C_A(r) = \sum_{k=0}^{n-1} G_k\,(A)\,\overline{G_{k+r}(A)} = \sum_{k=0}^{n-1} \left(\sum_{t=0}^{n-1} \Lambda_t\,(A)\,\omega^{kt}\right)\overline{\left(\sum_{s=0}^{n-1} \Lambda_s\,(A)\,\omega^{s(k+r)}\right)}$$

Reordering and using $\overline{\omega^{s(k+r)}} = \omega^{-s(k+r)}$:

$$C_A(r) = \sum_{t,s=0}^{n-1} \Lambda_t\,(A)\,\overline{\Lambda_s}(A)\,\omega^{-sr} \sum_{k=0}^{n-1} \omega^{k(t-s)}\,.$$

By orthogonality, $\sum_{k=0}^{n-1} \omega^{k(t-s)} = n\,\delta_{ts}$. Then:

$$C_A(r) = n \sum_{t=0}^{n-1} \Lambda_t\,(A)\,\overline{\Lambda_t(A)}\,\omega^{-tr} = n \sum_{t=0}^{n-1} |\,\Lambda_t(A)\,|^2\ \ \omega^{-tr}.$$

This is the DFT of the sequence $\{|\ \Lambda_t(A)\ |^2\}$ evaluated at $r$. Reindexing $s = t + r$ at $\mathbb{Z}_n$ yields the circular correlation form.

$$C_A(r) = n \sum_{t=0}^{n-1} \Lambda_t\,(A)\,\overline{\Lambda_{t+r}(A)}$$

## 9.3. Partial modal energy

**Theorem 9.7 (Partial modal energy)**

For $S \subseteq \{0, \dots, n-1\}$, let $E_S(A) := \sum_{k\in S} |\,G_k(A)\,|^2$. Then $E_S(A) \geq 0$ and $E_{\{0,\dots,n-1\}}(A) = \|\,D_\phi(A)\,\|_2^2$. If $A \in M_n(\mathbb{R})$ and $S$ is stable under $k \mapsto n-k$, then $E_S(A)$ is real and is calculated by grouping conjugate pairs:

$$E_S(A) = |\,G_0\,|^2 \cdot \mathbf{1}^1{}_{[0\in S]} + \sum_{k\in S, 1\leq k<n/2} 2\,|\,G_k(A)\,|^2 + |\,G_{n/2}\,|^2 \cdot \mathbf{1}_{[n/2\in S,\, n \text{ par}]}$$

Proof:

Each addend $|\ G_k(A)\ |^2$ is non-negative, therefore $E_S(A) \geq 0$. If $S \subseteq T \subseteq \mathbb{Z}_n$, then

---

[1] The symbol $\mathbf{1}_{[P]}$ denotes the indicator function of the proposition $P$: it holds 1 if $P$ it is true and 0 if $P$ it is false. For example, $\mathbf{1}_{[0\in S]} = 1$ when $0 \in S$and $\mathbf{1}_{[0\in S]} = 0$ when $0 \notin S$. This notation allows us to write in a single formula a sum whose extreme addends are present or absent depending on the conditions on $S$ and the parity of $n$, avoiding separation into cases.

$$E_T(A) - E_S(A) = \sum_{k \in T \setminus S} | G_k(A) |^2 \geq 0,$$

Therefore $E_S(A) \leq E_T(A)$. Taking $T = \mathbb{Z}_n$:

$$E_{\mathbb{Z}_n}(A) = \sum_{k=0}^{n-1} | G_k(A) |^2 = \| D_\phi(A) \|_2^2.$$

In the actual case, Proposition 3.8 gives $G_{n-k}(A) = \overline{G_k(A)}$, then $| G_{n-k}(A) |^2 = | G_k(A) |^2$. If $S$ is stable under $k \mapsto n-k$, the indices of $S \cap \{1, \ldots, \lfloor n/2 \rfloor\}$ are grouped into pairs $\{k, n-k\}$ that $2 \mid G_k \mid^2$ each contribute. The mode $k = 0$ (always fixed) contributes $| G_0 |^2$, and if $n$ is even, the mode $k = n/2$ (also fixed, and real by Proposition 3.8) contributes $| G_{n/2} |^2$. This gives exactly the formula of the statement, which is real and non-negative.

**9.4. Modal cofactor matrix**

**Definition 9.8 (Modal Cofactor Matrix)**

$$C^{(k)}(A) := (C_{ji}^{(k)}(A))_{i,j}.$$

For $k = 0$, $C^{(0)}(A)$ it coincides with the classical adjoint.

**Proposition 9.9.** For each fixed row $i$:

$$G_k(A) = \sum_j a_{ij} \; [C^{(k)}(A)]_{ji}$$

Proof:

By Theorem 5.8, for each fixed row $i$,

$$G_k(A) = \sum_{j=1}^{n} a_{ij} \; C_{ij}^{(k)}(A).$$

According to Definition 9.8, $C^{(k)}(A)_{ji} = C_{ij}^{(k)}(A)$. Therefore,

$$\sum_{j=1}^{n} a_{ij} \; [C^{(k)}(A)]_{ji} = \sum_{j=1}^{n} a_{ij} \; C_{ij}^{(k)}(A) = G_k(A).$$

The left-hand side is the dot product between the row $i$ and $A$ the $i$ column $C^{(k)}(A)$.

**9.5. Orbital spectral entropy**

**Definition 9.10 (Spectral distribution and entropy)**

Assuming $D_\varphi(A) \neq 0$, we define

$$p_k(A) := \frac{| G_k(A) |^2}{\| D_\varphi(A) \|_2^2}, k = 0, \dots, n-1,$$

so that $\sum_{k=0}^{n-1} p_k(A) = 1$. The orbital spectral entropy is defined by

$$H(A) := -\sum_{k=0}^{n-1} p_k(A) \log\ p_k(A),$$

with the convention $0 \log\ 0 := 0$.

**Example** ($n = 3$, $A$ tridiagonal). With $\| D_\varphi(A) \|_2^2 = 120$, the modal probabilities are

$$p_0 = \frac{16}{120} = \frac{2}{15} \approx 0.133, p_1 = p_2 = \frac{52}{120} = \frac{13}{30} \approx 0.433,$$

with $p_0 + p_1 + p_2 = 1$. Therefore,

$$H(A) = -(p_0 \log p_0 + 2\, p_1 \log p_1) \approx -(0.133 \log 0.133 + 2 \cdot 0.433 \log 0.433)$$

$\approx 0.993$ nat $\approx 90.4\%$ of $\log\ 3 \approx 1.099$.

**9.6. Vector Hadamard Inequality**

**Theorem 9.11 (Vector Hadamard Inequality).**

Let $A \in M_n(\mathbb{C})$ have columns $c_1, \dots, c_n$. Then

$$\| D_\phi(A) \|_2 \leq \sqrt{n!} \prod_{j=1}^{n} \| c_j \|_2.$$

Proof.

**Step 1.** By Corollary 9.3,

$$\| D_\phi(A) \|_2^2 = n \sum_{r=0}^{n-1} | \Lambda_r(A) |^2.$$

**Step 2.** Each class $\phi^{-1}(r)$ has $(n-1)!$ elements. By the Cauchy--Schwarz inequality,

$$| \Lambda_r(A) |^2 = | \sum_{\phi(\sigma)=r} m_\sigma(A) |^2 \leq (n-1)! \sum_{\phi(\sigma)=r} | m_\sigma(A) |^2.$$

**Step 3.** Summing over $r$, we obtain

$$\sum_{r=0}^{n-1} | \Lambda_r(A) |^2 \leq (n-1)! \sum_{\sigma \in S_n} | m_\sigma(A) |^2 = (n-1)!\ \mathrm{per}\ (| A |^{\circ 2}),$$

where $| A |^{\circ 2}$ denotes the matrix whose entries are $| a_{ij} |^2$.

For the structural and computational role of the permanent, see [9].

**Step 4.** We now bound the permanent:

$$\operatorname{per}(|A|^{\circ 2}) \le \prod_{j=1}^{n} \sum_{i=1}^{n} |a_{ij}|^2 = \prod_{j=1}^{n} \| c_j \|_2^2.$$

Indeed, by definition of the permanent,

$$\operatorname{per}(|A|^{\circ 2}) = \sum_{\sigma \in S_n} \prod_{j=1}^{n} |a_{\sigma^{-1}(j),j}|^2.$$

On the other hand, expanding the product

$$\prod_{j=1}^{n} \left( \sum_{i=1}^{n} |a_{ij}|^2 \right)$$

yields a sum over all maps

$\tau: \{1, \dots, n\} \to \{1, \dots, n\}$

of the terms

$$\prod_{j=1}^{n} |a_{\tau(j),j}|^2.$$

Since the permutations $S_n$ form a subset of all such maps, and all summands are nonnegative, it follows that

$$\operatorname{per}(|A|^{\circ 2}) \le \prod_{j=1}^{n} \sum_{i=1}^{n} |a_{ij}|^2 = \prod_{j=1}^{n} \| c_j \|_2^2.$$

**Step 5.** Combining the previous estimates, we obtain

$$\| D_\phi(A) \|_2^2 \le n \cdot (n-1)! \cdot \prod_{j=1}^{n} \| c_j \|_2^2 = n! \prod_{j=1}^{n} \| c_j \|_2^2.$$

Taking square roots gives

$$\| D_\phi(A) \|_2 \le \sqrt{n!} \prod_{j=1}^{n} \| c_j \|_2.$$

This proves the theorem.

**Corollary 9.11.1 (Comparison with the classical Hadamard inequality)**

The classical Hadamard inequality states that

$$| \det (A) | \leq \prod_{j=1}^{n} \| c_j \|_2,$$

where $c_j$ are the columns of $A$. Comparing with Theorem 9.11:

$$\| D_\phi(A) \|_2 \leq \sqrt{n!} \prod_{j=1}^{n} \| c_j \|_2 = \sqrt{n!} \cdot \text{(cota de Hadamard para } | \det A |).$$

Relationship between the two bounds. By Proposition 3.9, $| G_0(A) | = | \det (A) | \leq \mathrm{per}(| A |)$. By the Hadamard bound, $| \det A | \leq \prod_j \| c_j \|_2$. This does not directly imply Theorem 9.11, since $\| D_\phi \|_2$ it involves all modes, not just $G_0$. The proof of Theorem 9.11 requires going through orbital sums via Parseval and bounding by the permanent of $| A |^{\circ 2}$.

**The case of equality in vector Hadamard.** Inequality $\| D_\phi(A) \|_2 \leq \sqrt{n!} \prod_j \| c_j \|_2$ is equality only when:

1. There is equality in Cauchy-Schwarz within each class $\phi^{-1}(r)$: all monomials in that class have a constant modulus.
2. There is equality in the bound $\mathrm{per}(| A |^{\circ 2}) \leq \prod_j \| c_j \|_2^2$: each column concentrates its norm in a single entry (the column has a single non-zero element).

This corresponds, in practice, to the case of rescaled permutation matrices: $A = D \cdot P_\tau$ with $D$ diagonal. For such matrices, Corollary 3.15(b) gives $G_k(A) = \varepsilon(\tau) d_1 \cdots d_n \cdot \omega^{k\phi(\tau)}$ for all $k$, and $\| D_\phi(A) \|_2 = \sqrt{n} \cdot | \det A |$, with $\prod_j \| c_j \|_2 = \prod_j | d_{\tau^{-1}(j)} | = | \det A |$.

Proof of the Corollary: By Theorem 9.11, $\| D_\phi \|_2 \leq \sqrt{n!} \prod_j \| c_j \|_2$. The fundamental mode satisfies $| G_0 | = | \det A | \leq \| D_\phi \|_2$ (since $| G_0 |^2 \leq \sum_k | G_k |^2 = \| D_\phi \|_2^2$). Combining with the classical Hadamard bound $| \det A | \leq \prod_j \| c_j \|_2$: the vector bound is weaker by the factor $\sqrt{n!}$, but is more informative in that it bounds the entire norm of $D_\phi$, not just its fundamental component.

**9.7. Invariance under rescaling**

**Proposition 9.12.**

For all

$\lambda \in \mathbb{C}^*$,

$R_\infty(\lambda A) = R_\infty(A)$ and $V_2(\lambda A) = V_2(A)$,

provided that these quantities are defined.

**Proof.** By multilinearity (Theorem 3.6), each Leibniz term satisfies

$m_\sigma(\lambda A) = \lambda^n m_\sigma(A)$.

Consequently,

$\Lambda_r(\lambda A) = \lambda^n \Lambda_r(A), G_k(\lambda A) = \lambda^n G_k(A), D_\phi(\lambda A) = \lambda^n D_\phi(A)$.

Taking absolute values and norms, we obtain

$| G_0(\lambda A) | = | \lambda |^n | G_0(A) |, \| D_\phi(\lambda A) \|_2 = | \lambda |^n \| D_\phi(A) \|_2$,

and

$\text{per}\,(| \lambda A |) = | \lambda |^n \,\text{per}\,(| A |)$.

Therefore,

$$R_\infty(\lambda A) = \frac{| G_0(\lambda A) |}{\text{per}\,(| \lambda A |)} = \frac{| \lambda |^n | G_0(A) |}{| \lambda |^n \,\text{per}\,(| A |)} = R_\infty(A),$$

and

$$V_2(\lambda A) = 1 - \frac{| G_0(\lambda A) |^2}{\| D_\phi(\lambda A) \|_2^2} = 1 - \frac{| \lambda |^{2n} | G_0(A) |^2}{| \lambda |^{2n} \| D_\phi(A) \|_2^2} = V_2(A).$$

**As a consequence,**

$H(\lambda A) = H(A)$,

because

$p_k(\lambda A) = p_k(A)$ for every $k$.

### 9.8. Well-formulated conjectures

**Conjecture 9.13 (Gaussian concentration iid matrices)**

The following conjecture relies on the distributive symmetry of modes established in Corollary 8.4 and on analogies with classical concentration results in Gaussian determinants, but remains open: its proof would require controlling fourth-order moments and bounding the variance of $\| D_\phi(A) \|_2^2$.

For matrices $A \in M_n(\mathbb{R})$ with iid entries $\sim \mathcal{N}(0,1)$:

$$H(A) \underset{n\to\infty}{\overset{\text{prob.}}{\to}} \log n.$$

Motivation: Corollary 8.4 establishes the mean square value of the modes under the real iid Gaussian model. In particular, the distributive symmetry between the modal indices suggests a mean equipartition of the spectral energy. However, a concentration proof requires checking fourth-order and higher moments of the quantities $\mid G_k(A) \mid^2$, which leads to Wick–Isserlis-type Gaussian combinatorics and bounds on permanents of non-negative matrices.

**Conjecture 9.14 (Selective Modal Cancellation)**

The following conjecture describes a phenomenon of selective modal annulment observed in low-dimensional examples; its general proof, which would require solving an overdetermined polynomial system in the inputs of $A$, remains an open problem.

For $n \geq 3$, given $c \in \mathbb{C}$ and $T \subseteq \{1, \dots, n-1\}$, there exist nontrivial families of matrices $A \in M_n(\mathbb{C})$ with

$G_0(A) = c$ y $G_k(A) = 0 \quad \forall k \in T.$

The existence of such families would generalize the classical problem $\det(A) = c$ to an orbital vector system of simultaneous polynomial equations.

## 10. Finite Trigonometric Interpolation and Analytic Extension

### 10.1. Interpolated function

**Definition 10.1 (Interpolated function of the determinant vector)**

The vector $D_\phi(A)$ provides $n$ discrete values $G_k(A)$ indexed by $\mathbb{Z}_n$. A natural question is whether there exists an analytic function that interpolates these values and allows evaluation at non-integer frequencies $\xi \in \mathbb{R}$. The answer is yes: the finite trigonometric function $G(z; A) = \sum_r \Lambda_r(A) e^{2\pi i r z/n}$ extends $D_\phi(A)$ to the complex plane with all the properties developed in this section.

$$G(z;\, A) := \sum_{r=0}^{n-1} \Lambda_r(A)\, e^{2\pi i r z/n}, z \in \mathbb{C}. \qquad (6)$$

This is the natural finite trigonometric interpolation of $D_\phi(A)$: a linear combination of $n$ exponentials with coefficients $\Lambda_r(A)$.

**Theorem 10.2 (Basic properties of interpolation).**

The function $G(\cdot; A)$ satisfies:

- **Exact interpolation:**

$G(k; A) = G_k(A)$ for $k = 0,1, \dots, n-1$.

- **Periodicity:**

$G(z + n; A) = G(z; A)$ for all $z \in \mathbb{C}$.

- **Entirety:** $G(\cdot; A)$ is an entire function on $\mathbb{C}$.

- **Exponential growth:**

$$| G(z; A) | \le C_A \, e^{2\pi(n-1)|\mathrm{Im}\ z|/n}, C_A := \sum_{r=0}^{n-1} | \Lambda_r(A) |.$$

**Proof.**

**(1)** By definition,

$$G(k; A) = \sum_{r=0}^{n-1} \Lambda_r \, (A) e^{2\pi i r k/n} = \sum_{r=0}^{n-1} \Lambda_r \, (A) \omega^{kr} = G_k(A).$$

**(2)** We have

$$e^{2\pi i r(z+n)/n} = e^{2\pi i r z/n} \, e^{2\pi i r} = e^{2\pi i r z/n},$$

since

$$e^{2\pi i r} = 1.$$

Hence

$$G(z + n; A) = G(z; A).$$

**(3)** The function $G(\cdot; A)$ is a finite sum of entire functions $z \mapsto e^{2\pi i r z/n}$, and is therefore entire.

**(4)** Let $z = x + iy$. Then

$$| e^{2\pi i r z/n} | = e^{-2\pi r y/n}.$$

If $y > 0$, the maximum over $r \in \{0, \dots, n-1\}$ is attained at $r = 0$, giving the value 1. If $y < 0$, the maximum is attained at $r = n - 1$, giving

$$e^{2\pi(n-1)|y|/n}.$$

Therefore,

$$| G(z; A) | \le \sum_{r=0}^{n-1} | \Lambda_r(A) | \quad | e^{2\pi i r z/n} | \le \left( \sum_{r=0}^{n-1} | \Lambda_r(A) | \right) e^{2\pi(n-1)|\mathrm{Im}\ z|/n}.$$

Thus

$$| G(z; A) | \le C_A \, e^{2\pi(n-1)|\mathrm{Im}\ z|/n}.$$

### 10.2. Interpolated value of continuous frequency

Definition 10.3 (Interpolated continuous frequency value)

For $\xi \in \mathbb{R}$:

$$G_\xi(A) := G(\xi;\ A) = \sum_{r=0}^{n-1} \Lambda_r\,(A)\, e^{2\pi i r\xi/n}. \qquad (7)$$

This object extends discrete modes $\{G_k(A)\}_{k=0}^{n-1}$to a continuous parameter $\xi \in \mathbb{R}$. For $\xi = p/q \in \mathbb{Q}$, the exponentials $e^{2\pi i r p/(nq)}$ are roots of unity of order $nq$.

**Proposition 10.4 (Properties of the interpolated value)**

For all $\xi \in \mathbb{R}$:

1. $G_{\xi+n}(A) = G_\xi(A)$ (periodicity).
2. $G(k;\ A) = G_k(A)$ for $k \in \{0, \dots, n-1\}$ integer (the interpolated continuous function recovers exactly the discrete modes).
3. If $n$ even: $G_{n/2}(A) = \sum_{r=0}^{n-1} \Lambda_r\,(A) \cdot (-1)^r$ (Nyquist value).

Proof:

**(1)** By the periodicity of Theorem 10.2: $G_{\xi+n}(A) = G(\xi + n; A) = G(\xi; A) = G_\xi(A)$.

**(2)** For $k \in \{0, \dots, n-1\}$ integer, $G(k; A) = \sum_r \Lambda_r\,(A)\, e^{2\pi i r k/n} = \sum_r \Lambda_r\,(A)\,\omega^{kr} = G_k(A)$, which is exactly property (1) of Theorem 10.2.

**(3)** Evaluating at $\xi = n/2$: $e^{2\pi i r(n/2)/n} = e^{\pi i r} = (\cos\,\pi r + i\sin\,\pi r) = (-1)^r$, then $G_{n/2}(A) = \sum_{r=0}^{n-1} \Lambda_r\,(A)\,(-1)^r$.

### 10.3. Symmetry for real matrices

**Proposition 10.5 (Continuous Conjugate Symmetry)**

If $A \in M_n(\mathbb{R})$, then

$\overline{G(\xi;\ A)} = G(-\xi;\ A) = G(n-\xi;\ A)$ For all. $\xi \in \mathbb{R}$

In particular, $G(0; A) = \det\,(A) \in \mathbb{R}$ and, if $n$ it is even, $G(n/2; A) \in \mathbb{R}$.

Proof: For $A$ real, $\Lambda_r(A) \in \mathbb{R}$. Then:

$$\overline{G(\xi; A)} = \sum_r \Lambda_r\,(A)\, e^{-2\pi i r\xi/n} = G(-\xi; A).$$

Equality $G(-\xi; A) = G(n-\xi; A)$ follows from periodicity $G(-\xi; A) = G(-\xi+n; A) = G(n-\xi; A)$.

This proposition is the continuous version of the discrete symmetry $G_{n-k}(A) = \overline{G_k(A)}$ of Proposition 3.8.

**10.4. Derivatives and orbital moments**

**Theorem 10.6 (Derivatives of interpolation)**

For all $p \geq 0$:

$$\frac{d^p}{dz^p} G(z; A) = \sum_{r=0}^{n-1} \Lambda_r(A) \left(\frac{2\pi i r}{n}\right)^p e^{2\pi i r z/n}.$$

In $z = 0$:

$$G'(0; A) = \frac{2\pi i}{n} \sum_{r=0}^{n-1} r\ \Lambda_r(A), G''(0; A) = -\frac{4\pi^2}{n^2} \sum_{r=0}^{n-1} r^2\ \Lambda_r(A).$$

Proof: Term-by-term differentiation of the finite sum (6) is valid (entire functions).

The function $G(z; A) = \sum_{r=0}^{n-1} \Lambda_r(A)\, e^{2\pi i r z/n}$ is a finite sum of entire functions. For finite sums of differentiable functions, differentiation can be done term by term:

$$\frac{d^p}{dz^p} e^{2\pi i r z/n} = \left(\frac{2\pi i r}{n}\right)^p e^{2\pi i r z/n}.$$

Adding in $r$gives the general formula. Evaluating at $z = 0$ and using $e^0 = 1$ gives the expressions for $G'(0; A)$ and $G''(0; A)$.

**Interpretation.** The derivatives of $G$ in $z = 0$ are the discrete moments of the orbital distribution $\{\Lambda_r(A)\}$. The first derivative measures the orbital "center of mass"; the second, the orbital "variance".

**Definition 10.7 (Continuous moments)**

For $p \geq 0$:

$$\mu_p(A) := \frac{1}{n} \int_0^n |\ G(\xi; A)\ |^2 \cdot \xi^p\, d\xi.$$

They quantify the distribution of the interpolated spectral energy throughout the cycle $[0, n)$.

**10.5. Continuous Parseval Identity**

**Theorem 10.8 (Continuous Parseval)**

We have the chain of equalities:

$$\frac{1}{n}\int_0^n \mid G(\xi;A) \mid^2 \ d\xi = \sum_{r=0}^{n-1} \mid \Lambda_r(A) \mid^2 = \frac{1}{n} \parallel D_\phi(A) \parallel_2^2.$$

Proof: We expand:

$$\mid G(\xi;A) \mid^2 = \sum_{r,s=0}^{n-1} \Lambda_r\,(A)\,\bar{\Lambda_s}(A)\,e^{2\pi i(r-s)\xi/n}.$$

Integrating term by term in $[0, n]$:

$$\frac{1}{n}\int_0^n e^{2\pi i(r-s)\xi/n}\ d\xi = \begin{cases} 1 & r = s, \\ 0 & r \neq s. \end{cases}$$

Therefore: $\frac{1}{n}\int_0^n \mid G(\xi;A) \mid^2 \ d\xi = \sum_r \mid \Lambda_r(A) \mid^2$.

The second equality is Corollary 9.3: $\sum_r \mid \Lambda_r \mid^2 = \frac{1}{n} \parallel D_\phi \parallel_2^2$.

**Remark 10.9** (Three unified perspectives). Theorem 10.8 unites:

- **(a)** The continuous integral of the energy of the integer function $G(\cdot;A)$.
- **(b)** The sum of squares of the orbital sums $\{\Lambda_r(A)\}$.
- **(c)** The Euclidean norm of the determinant vector $D_\phi(A)$.

Verification: $(1/3) \cdot \int_0^3 \mid G(\xi;A) \mid^2 d\xi = 40 = \parallel D_\phi \parallel_2^2/3 = 120/3$.

**10.6. The orbital polynomial**

Definition 10.10 (Orbital polynomial).

$$P_A(q) := \sum_{r=0}^{n-1} \Lambda_r\,(A)\,q^r \ \in \ \mathbb{C}[q].$$

The interpolated function satisfies

$$G(z;A) = P_A\left(e^{2\pi iz/n}\right).$$

That is, $G(\cdot;A)$ it is the composition of the polynomial $P_A$ with the exponential $e^{2\pi iz/n}$.

**Proposition 10.11 (Zeros of interpolation)**

We denote by

$G(z;A) = 0 \ \Leftrightarrow \ e^{2\pi iz/n}$ es raíz de $P_A$.

In particular, if $G(\cdot;A) \not\equiv 0$:

- The real zeros of $G$in $[0, n)$ are finite (by analyticity).
- They correspond exactly to the roots of $P_A$ on the unit circle $| q |= 1$.

Proof:

By Definition 10.10, $G(z; A) = P_A(e^{2\pi iz/n})$. Therefore:

$G(z; A) = 0 \iff P_A\left(e^{2\pi iz/n}\right) = 0 \iff e^{2\pi iz/n}$ es raíz de $P_A$.
This establishes the equivalence For all $z \in \mathbb{C}$.

**Real case** $z = x \in \mathbb{R}$**:** Since $| e^{2\pi ix/n} |= 1$, the real zeros of $G(\cdot; A)$ correspond exactly to the roots of $P_A$ on the unit circle $\mathbb{S}^1 = \{| q |= 1\}$.

**Finitude:** $P_A(q) = \sum_{r=0}^{n-1} \Lambda_r (A) q^r$ is a polynomial of degree at most $n - 1$. If $P_A \equiv/0$, has at most $n - 1$ roots in $\mathbb{C}$, then at most $n - 1$ roots on $\mathbb{S}^1$. Since $e^{2\pi ix/n}$ is periodic modulo $n$, $x$ it suffices to study $x \in [0, n)$, where the map $x \mapsto e^{2\pi ix/n}$ traverses $\mathbb{S}^1$ exactly once. This yields a finite number of real zeros in $[0, n)$

**Observation 10.12.** Perfect cancellation $G_\xi(A) = 0$ for some $\xi \in \mathbb{R}$ (with $G_0(A) = \det(A) \neq 0$) defines a notion of spectral singularity: the matrix is not singular in the classical sense, but possesses an orbital phase configuration that produces perfect cancellation at a non-integer frequency. This opens a new dimension of matrix classification by "distance to singularity" in phase, not just in magnitude.

**10.7. Examples**

**Case** $n = 2$**:** $A = \begin{pmatrix} a & b \\ c & d \end{pmatrix}$.

$\Lambda_0 = ad$, $\Lambda_1 = -bc$. Interpolated function:

$G(z; A) = ad - bc \cdot e^{\pi iz}$.
Integer values: $G(0; A) = ad - bc = \det(A)$, $G(1; A) = ad + bc = \mathrm{per}(A)$.

Fractional value: $G(1/2; A) = ad - bc \cdot e^{\pi i/2} = ad - i\, bc$.

The value $G_{1/2}(A) = ad - i\, bc$ combines the main diagonal with the antidiagonal through a complex phase $e^{i\pi/2} = i$, producing an invariant algebraically distinct from both the determinant and the mode $G_1$.

Orbital polynomial: $P_A(q) = ad - bc \cdot q$. Its root is $q_0 = ad/(bc)$(when $bc \neq 0$). The determinant vanishes if and only if $q_0 = 1$if and only if $ad = bc$; a real zero of $G$ exists if and only if $\mid q_0 \mid = 1$ if and only if $\mid ad \mid = \mid bc \mid$.

**Case $n = 3$:$A$ tridiagonal.**

$\Lambda_0 = 6, \Lambda_1 = -2, \Lambda_2 = 0$. Interpolated function:

$G(z; A) = 6 - 2\, e^{2\pi i z/3}$.

**Verification in integers:**

$$G(0; A) = 6 - 2 = 4\ , G(1; A) = 6 - 2\omega = 7 - i\sqrt{3}\ , G(2; A) = 6 - 2\omega^2 = 7 + i\sqrt{3}\ .$$

**Fractional value** in $\xi = 3/2$:

$G(3/2; A) = 6 - 2\, e^{2\pi i \cdot (3/2)/3} = 6 - 2\, e^{\pi i} = 6 - 2 \cdot (-1) = 8.$

This value $G_{3/2}(A) = 8$ does not coincide with $\det(A) = 4$ either $\mid G_1(A) \mid = \sqrt{52} \approx 7.21$: continuous interpolation detects invisible orbital combinations in discrete modes.

**Orbital polynomial:** $P_A(q) = 6 - 2q$ its only root is $q_0 = 3$, with $\mid q_0 \mid = 3 > 1$. Then $G(\cdot\, ; A)$ it has no real zeros in $[0,3)$, consistent with the non-singularity of $A$($\det(A) = 4 \neq 0$).

### 10.8. Open problems

**Problem 10.13 (Inverse Spectral Synthesis)**

Given $f(z) = \sum_{r=0}^{n-1} c_r\, e^{2\pi i r z/n}$ arbitrary $c_r \in \mathbb{C}$ numbers, do there exist numbers $A \in M_n(\mathbb{C})$ and $\phi$ such that $G(z; A) = f(z)$? Equivalently: $A \mapsto \{\Lambda_r(A)\}_{r=0}^{n-1}$ is the map surjective onto $\mathbb{C}^n$?

**Problem 10.14 (Zeros on the unit circle)**

Characterize, in terms of the matrix structure of $A$, when $P_A(q)$ it has roots over $\mid q \mid = 1$. What does the existence of $\xi \in \mathbb{R}$ with $G(\xi; A) = 0$ over the distribution of Leibniz terms in the orbits imply?

**Problem 10.15 (Stability of zeros under disturbances)**

For $A \mapsto A + \varepsilon H$ fixed $H \in M_n(\mathbb{C})$, describe the movement of the roots of $P_{A+\varepsilon H}$ as a function of $\varepsilon \in \mathbb{C}$. Are the real zeros of $G(\cdot; A)$ structurally stable under small matrix perturbations?

## 11. Canonicity of the Orbital Phase, Covariance of the Vector Determinant and Structural Limit of the Circulating Case

### 11.1. Orbital Decomposition and Uniqueness of the Canonical Representative

The results in the previous sections depend on the choice of the phase function $\phi$. This section answers precisely: to what extent is it $\phi$ unique?, how does it transform $D_{\phi}(A)$ when $\phi$ it changes?, and what quantities are genuinely independent of that choice?

Let be $S_n$ the symmetric group and let

$\rho = (1\ 2\ \cdots\ n) \in S_n$

the standard cycle. We consider the rightward action $C_n = \langle \rho \rangle$ on $S_n$:

$(\sigma, r) \longmapsto \sigma \circ \rho^r.$

According to Proposition 6.1 of the manuscript, this action is free and partitions $S_n$ into $(n - 1)!$orbits of size $n$.

**Definition 11.1 (Cyclic Orbit).** The orbit $\sigma \in S_n$ under the right action of $C_n$ is the set

$\mathcal{O}(\sigma) = \{\sigma \circ \rho^r : r \in \mathbb{Z}_n\}.$

The first question of canonicity consists of justifying that the anchoring $\sigma^*(1) = 1$ selects a single representative per orbit.

**Proposition 11.2 (Uniqueness of the canonical representative)**

Each orbit $\mathcal{O} \subset S_n$ contains a unique element $\sigma^*$ such that $\sigma^*(1) = 1$.

Proof: Let $\sigma \in \mathcal{O}$. We want to find $r \in \mathbb{Z}_n$ such that $(\sigma \circ \rho^r)(1) = 1$, that is,

$\sigma(\rho^r(1)) = 1.$

Since $\rho^r(1) = r + 1(\text{mod} n)$ (with the convention $\{1, \dots, n\}$), the condition becomes $\sigma(r + 1) = 1$. Because $\sigma$ it is bijective, there exists a unique index $j \in \{1, \dots, n\}$ with $\sigma(j) = 1$, namely $j = \sigma^{-1}(1)$. Then there exists a unique $r \in \mathbb{Z}_n$ such that $r + 1 \equiv \sigma^{-1}(1)\ (\text{mod} n)$, that is,

$r \equiv \sigma^{-1}(1) - 1(\text{mod} n).$

That value $r$ produces the only element $\sigma^* = \sigma \circ \rho^r$ with $\sigma^*(1) = 1$.

**Uniqueness.** If $\sigma \circ \rho^r$and $\sigma \circ \rho^s$both satisfy the condition $(\cdot)(1) = 1$, then $\sigma(\rho^r(1)) = \sigma(\rho^s(1)) = 1$. Since $\sigma$ is injective, $\rho^r(1) = \rho^s(1)$, then $r \equiv s(\text{mod} n)$.

**Corollary 11.3 (Algorithmic formula of the representative)**

Given $\sigma \in S_n$, the canonical representative of its orbit is

$\sigma^* = \sigma \circ \rho^{r^*}, r^* \equiv \sigma^{-1}(1) - 1(\text{mod} n).$

Its orbital position is $\phi(\sigma) = r^*$

**11.2. Canonical Phase and Uniqueness of Orbital Parametrization**

**Definition 11.4 (Phase function)**

A phase function is a map $\phi: S_n \to \mathbb{Z}_n$ such that, for each orbit $\mathcal{O}(\sigma^*)$ with canonical representative $\sigma^*(\ \sigma^*(1) = 1)$, we have

$\phi(\sigma^*) = 0, \phi(\sigma^* \circ \rho^r) = r \forall\, r \in \mathbb{Z}_n$.

**Proposition 11.5 (Uniqueness of the canonical phase)**

Once the anchor $\sigma^*(1) = 1$ is fixed and secured $\rho$, there is a single function $\phi_{\text{can}}: S_n \to \mathbb{Z}_n$ that satisfies Definition 11.4.

Proof: Let $\sigma \in S_n$. By Proposition 11.2, the orbit $\mathcal{O}(\sigma)$ contains a unique canonical representative $\sigma^*$ with $\sigma^*(1) = 1$. Since the action of $C_n$ is free (Proposition 6.1), there exists a unique $r \in \{0, \dots, n-1\}$ such that $\sigma = \sigma^* \circ \rho^r$. We then define

$\phi_{\text{can}}(\sigma) := r$.

It satisfies Definition 11.4: if $\sigma = \sigma^*$ (the canonical representative), then $r = 0$ and $\phi_{\text{can}}(\sigma^*) = 0$. If $\sigma = \sigma^* \circ \rho^r$, then $\phi_{\text{can}}(\sigma) = r$ by construction.

Uniqueness: Let be $\psi: S_n \to \mathbb{Z}_n$ another function that satisfies Definition 11.4. For any $\sigma = \sigma^* \circ \rho^r$, the property requires $\psi(\sigma^* \circ \rho^r) = r$. Since this decomposition is unique, $\psi(\sigma) = r = \phi_{\text{can}}(\sigma)$ for all $\sigma$. Then $\psi = \phi_{\text{can}}$.

**11.3. Two Levels of Phase Change: Global Gauge and Orbital Gauge**

Although the canonical phase is unique once the anchor is fixed, it is convenient to separate two types of phase transformations, as they have very different effects on $D_\phi(A)$.

**11.3.1. Global Gauge Affine**

**Definition 11.6.**

For $u \in (\mathbb{Z}_n)^\times$ (invertible integers modulo $n$) and $\beta \in \mathbb{Z}_n$, the global affine gauge is the transformation

$\phi_{u,\beta}(\sigma) := u\,\phi(\sigma) + \beta (\text{mod} n)$.

**11.3.2. Orbital Gauge per Orbit**

**Definition 11.7.**

Let be $\alpha$: {órbitas de $S_n$} $\to \mathbb{Z}_n$ a function. The associated orbital gauge is

$\phi_\alpha(\sigma) := \phi(\sigma) + \alpha(\mathcal{O}(\sigma))(\mathrm{mod} n)$.

This change rotates the phases of each orbit separately, with possible different rotations for each one.

**11.4. Transformation Law under Global Gauge Affine**

Recall that

$$G_k^{(\phi)}(A) = \sum_{\sigma \in S_n} m_\sigma\,(A)\,\omega^{k\phi(\sigma)}, D_\phi(A) = (G_0^{(\phi)}(A), \dots, G_{n-1}^{(\phi)}(A)).$$

**Theorem 11.8 (Covariance under global affine gauge).**

Let

$\phi' = \phi_{u,\beta} := u\phi + \beta, u \in (\mathbb{Z}/n\mathbb{Z})^\times, \beta \in \mathbb{Z}/n\mathbb{Z}.$

Then, for all $k \in \mathbb{Z}/n\mathbb{Z}$,

$G_k^{\phi'}(A) = \omega^{k\beta}\, G_{uk}^{\phi}(A).$

In matrix form,

$D_{\phi'}(A) = U_\beta\, P_u\, D_\phi(A),$

where

$U_\beta = \mathrm{diag}\,(1, \omega^\beta, \omega^{2\beta}, \dots, \omega^{(n-1)\beta}),$

and $P_u$ is the permutation matrix induced by

$k \mapsto uk(\mathrm{mod} n).$

**Proof.** By definition,

$$G_k^{\phi'}(A) = \sum_\sigma m_\sigma\,(A)\,\omega^{k\phi'(\sigma)} = \sum_\sigma m_\sigma\,(A)\,\omega^{k(u\phi(\sigma)+\beta)}.$$

Factoring out the constant term $\omega^{k\beta}$, we obtain

$$G_k^{\phi'}(A) = \omega^{k\beta} \sum_\sigma m_\sigma\,(A)\,\omega^{ku\phi(\sigma)} = \omega^{k\beta}\, G_{uk}^{\phi}(A).$$

The matrix formula is the coordinatewise rewriting of this identity: the $k$-th component of $D_{\phi'}(A)$ is $\omega^{k\beta}$ times the $uk$-th component of $D_\phi(A)$, which is exactly the action of $U_\beta P_u$.

**11.5. Invariants and Covariants under Global Gauge**

**Corollary 11.9 (Invariants under affine global gauge)**

Under the change $\phi \mapsto \phi_{u,\beta}$, the following are invariant:

1. $G_0(A) = \det(A)$.
2. $\| D_\phi(A) \|_2$.
3. The multiset of moduli $\{| G_0(A) |, \dots, | G_{n-1}(A) |\}$.
4. The orbital spectral entropy $H(A) = -\sum_k p_k \log p_k$, $p_k = | G_k |^2 / \| D_\phi \|_2^2$.
5. The index $V_2(A) = 1 - | G_0(A) |^2 / \| D_\phi(A) \|_2^2$.

Proof: By Theorem 11.8, $G_k^{(\phi')}(A) = \omega^{k\beta} G_{uk}^{(\phi)}(A)$, so that $| G_k^{(\phi')}(A) | = | G_{uk}^{(\phi)}(A) |$.

**Item 1:** for $k = 0$: $G_0^{(\phi')}(A) = \omega^0 G_0^{(\phi)}(A) = \det(A)$.

**Item 2:** From the above, $\| D_{\phi'}(A) \|_2^2 = \sum_k | G_k^{(\phi')} |^2 = \sum_k | G_{uk}^{(\phi)} |^2 = \sum_k | G_k^{(\phi)} |^2 = \| D_\phi(A) \|_2^2$.

**Item 3** : the map $k \mapsto uk \ mod n$ is a bijection of $\mathbb{Z}_n$ onto itself (since $u \in (\mathbb{Z}_n)^\times$). Therefore, the set of moduli $\left\{| G_0^{(\phi')} |, \dots, | G_{n-1}^{(\phi')} |\right\}$ is a permutation of the set $\left\{| G_0^{(\phi)} |, \dots, | G_{n-1}^{(\phi)} |\right\}$.

**Item 4:** The weights $p_k^{(\phi')}(A) = | G_k^{(\phi')} |^2 / \| D_{\phi'} \|_2^2 = | G_{uk}^{(\phi)} |^2 / \| D_\phi \|_2^2 = p_{uk}^{(\phi)}(A)$ are a permutation of the original weights. The entropy $H = -\sum_k p_k \log p_k$ is symmetric in its arguments, therefore $H^{(\phi')}(A) = H^{(\phi)}(A)$.

**Item 5:** $V_2^{(\phi')}(A) = 1 - | G_0^{(\phi')} |^2 / \| D_{\phi'} \|_2^2 = 1 - | \det(A) |^2 / \| D_\phi \|_2^2 = V_2^{(\phi)}(A)$.

### 11.6. Concrete example ( $n = 3$): global affine gauge and orbital gauge

Let $A = \begin{pmatrix} 2 & 1 & 0 \\ 1 & 2 & 1 \\ 0 & 1 & 2 \end{pmatrix}$ be the tridiagonal of §4.2: Under the canonical phase $\phi$:

$$D_\phi(A) = \left(4,\ 7 - i\sqrt{3},\ 7 + i\sqrt{3}\right), \qquad |D_\phi(A)|_2^2 = 120$$

**Global affine gauge** $u = 2, \beta = 1$(with $n = 3$, $u = 2 \in (\mathbb{Z}_3)^\times$, $\beta = 1$):

By Theorem 11.8, $G_k^{(\phi')}(A) = \omega^k G_{2k}^{(\phi)}(A)$. Calculating:

$$G_0^{(\phi')} = \omega^0 G_0^{(\phi)} = 4 = \det(A)$$

$G_1^{(\phi')} = \omega^1\, G_2^{(\phi)} = \omega(7 + i\sqrt{3}), G_2^{(\phi')} = \omega^2\, G_4^{(\phi)} = \omega^2\, G_1^{(\phi)} = \omega^2(7 - i\sqrt{3})$.

The moduli: $\mid G_1^{(\phi')} \mid = \mid 7 + i\sqrt{3} \mid = \sqrt{52}$, $\mid G_2^{(\phi')} \mid = \sqrt{52}$, equal to the originals. Invariants preserved: $\| D_{\phi'} \|_2^2 = 120$, $H(A)$, $V_2(A)$, verifying Corollary 11.9.

**Now consider a non-uniform orbital gauge** $c(\mathcal{O}_1) = 1, c(\mathcal{O}_2) = 0$ (there $n = 3$ are 2 orbits):

As stated in Observation 11.11 below, the modes are affected unevenly according to the contribution of each orbit. The mode $G_1$ receives weight $\omega^{k\cdot 1}$ from $\mathcal{O}_1$ and weight $\omega^{k\cdot 0}$ from $\mathcal{O}_2$, producing a $\tilde{G}_1 \neq G_1$ general. In contrast, $G_0 = \det(A) = 4$ remains invariant. This illustrates that outside the global affine gauge, only $G_0$ is intrinsic.

### 11.7. Genuine Dependence under Orbital Gauge

**Definition 11.10.** For each orbit $\mathcal{O}$ and each $k$, we define the orbital contribution to the mode $k$:

$$H_k(\mathcal{O}; A) := \sum_{\sigma \in \mathcal{O}} m_\sigma(A)\, \omega^{k\phi(\sigma)}.$$

So $G_k^{(\phi)}(A) = \sum_{\mathcal{O}} H_k(\mathcal{O}; A)$.

**Observation 11.11 (Dependence under orbital gauge).**

Let

$c: \{\text{orbits}\} \to \mathbb{Z}/n\mathbb{Z}$

be a function on the set of orbits, and define a new phase function by

$\tilde{\phi}(\sigma) := \phi(\sigma) + c(\mathcal{O}_\sigma)$,

where $\mathcal{O}_\sigma$ denotes the orbit of $\sigma$. Then

$$\tilde{G}_k(A) = \sum_{\mathcal{O}} \omega^{k\, c(\mathcal{O})}\; G_k^{\mathcal{O}}(A),$$

where $G_k^{\mathcal{O}}(A)$ denotes the contribution of the orbit $\mathcal{O}$ to the mode $G_k(A)$.

Consequently, unlike the global affine gauge, an orbital gauge does not generally act by a single global permutation of the modes together with a common phase factor. Its effect depends on the actual orbital decomposition of $A$. Therefore, quantities such as the modal distribution $V_2(A)$ and the spectral entropy $H(A)$ may vary when the orbital gauge changes.

**Comment.** This observation does not claim a complete classification of all possible transformations under orbital gauge. It only makes clear that, outside the global affine case, the dependence on the phase convention is genuine and does not generally reduce to a uniform covariance independent of $A$.

**Conclusion.** The theory has two levels:

- Under global affine gauge, the formalism is covariant and its spectral invariants are well defined.
- Under arbitrary orbital gauge, only

$G_0(A) = \det(A)$
remains universally intrinsic.

**11.8. The structural obstacle in the circulant case: Incompatibility of Degrees**

**Theorem 11.12 (Structural impossibility of identity $G_k(C) = \lambda_k(C)$)**

Let $n \geq 2$. Given any phase function $\phi$ and any $k \in \{0, \dots, n-1\}$, the function

$c = (c_0, \dots, c_{n-1}) \longmapsto G_k(\mathrm{circ}(c_0, \dots, c_{n-1}))$
is a homogeneous polynomial of degree $n$ in the variables $c_j$. In contrast,

$$c \longmapsto \lambda_k(\mathrm{circ}(c)) = \sum_{j=0}^{n-1} c_j \ \omega^{kj}$$

is homogeneous of degree 1. Therefore, the identity $G_k(C) = \lambda_k(C)$ for all circulant $C$ is impossible for $n \geq 2$.

Proof: Each Leibniz term

$$m_\sigma(A) = \varepsilon(\sigma) \prod_{i=1}^{n} a_{i,\sigma(i)}$$

is homogeneous of degree $n$in the entries of $A$. Then each orbital sum $\Lambda_r(A)$ and each mode

$$G_k(A) = \sum_{r=0}^{n-1} \Lambda_r\,(A)\,\omega^{kr}$$

It is also homogeneous of degree $n$. Restricting to $C = \mathrm{circ}(c_0, \dots, c_{n-1})$, whose entries are linear in the $c_j$, it turns out that $G_k(C)$ it is a homogeneous polynomial of degree $n$ in the $c_j$.

On the other hand, $\lambda_k(C) = \sum_j c_j\,\omega^{kj}$ it is linear (homogeneous of degree 1) in the $c_j$. Two homogeneous polynomials of different degrees cannot coincide identically in $\mathbb{C}^n$ unless both are zero, which does not occur here.

**11.9. What Should "Additional Standardization" Mean?**

**Corollary 11.13.**

If one seeks an identity of the type

$\mathcal{N}_k(C)\, G_k(C) = \lambda_{\pi(k)}(C)$

valid for all circulant currency $C$ (with $\pi$ a fixed permutation of the indices), then $\mathcal{N}_k(C)$ it must be homogeneous of degree $1-n$.

Proof: Scaling by $t \in \mathbb{C}^*$: $G_k(tC) = t^n G_k(C)$ and $\lambda_{\pi(k)}(tC) = t\, \lambda_{\pi(k)}(C)$. Substituting:

$\mathcal{N}_k(tC)\, t^n G_k(C) = t\, \lambda_{\pi(k)}(C) \implies \mathcal{N}_k(tC) = t^{1-n} \mathcal{N}_k(C).$

**Interpretation.** The desired normalization cannot be a simple constant factor: it must transform a homogeneous quantity of degree $n$ into one of degree 1, and therefore must scale with total degree $1-n$. In particular, any exact identity of this type will require a nonlinear normalization dependent on the matrix. The precise analysis of its regularity and its behavior near singular matrices remains open.

**11.10. What Is Actually True in the Circulant Case**

**Proposition 11.14 (Character Compatibility)**

Let $\widehat{C_n}$ the dual group of $C_n$. Both modes

$$G_k(A) = \sum_{r=0}^{n-1} \Lambda_r\,(A)\, \omega^{kr}$$

such as the spectral diagonalization of a circulating

$$\lambda_k(C) = \sum_{j=0}^{n-1} c_j\ \omega^{kj}$$

They are indexed by the same character system $\chi_k(r) = \omega^{kr}$, $k \in \mathbb{Z}_n$.

Proof: By Definition 3.2, $G_k(A) = \sum_{r=0}^{n-1} \Lambda_r\,(A)\, \omega^{kr}$, which is the DFT of the sequence $\{\Lambda_r(A)\}_{r=0}^{n-1}$evaluated at frequency $k$, using the characters $\chi_k(r) = \omega^{kr}$of $\mathbb{Z}_n$. On the other hand, the classical diagonalization of circulants gives $\lambda_k(C) = \sum_{j=0}^{n-1} c_j\ \omega^{kj}$, which is the DFT of the first row $(c_0, \dots, c_{n-1})$ evaluated at $k$, with the same characters $\chi_k$. The structural difference is that the coefficients are different: orbital sums $\Lambda_r(A)$ in one case, entries from the first row $c_j$ in the other.

**Observation 11.14.1 (Why the same DFT does not imply the same objects)**

Proposition 11.14 states that both $G_k(A)$ are $\lambda_k(C)$ DFTs on $\mathbb{Z}_n$ the same characters $\chi_k(r) = \omega^{kr}$. But their inputs are radically different in nature:

- **Inputs of $\lambda_k(C)$:** the coefficients $c_j$ of the first row of the circulating, which are linear in the inputs of $C$ (degree 1).
- **Entries of $G_k(A)$:** the orbital sums $\Lambda_r(A)$, which are polynomials of degree $n$ in the entries of $A$(Theorem 3.6 and Theorem 11.12).

This difference in degree is algebraically insurmountable: if $F: \mathbb{C}^{n\times n} \to \mathbb{C}$ is a polynomial of degree $d$and $G: \mathbb{C}^{n\times n} \to \mathbb{C}$ is a polynomial of degree $d'$ with $d \neq d'$, then $F(A) = G(A)$ for all $A \in \mathbb{C}^{n\times n}$ would imply $F - G \equiv 0$, which is impossible for two homogeneous polynomials of different degrees over a positive dimension space.

**Analogy with discrete Fourier.** Imagine we have the same orthonormal Fourier basis $\{e^{2\pi ikr/n}\}$ to expand two different functions $f, g: \mathbb{Z}_n \to \mathbb{C}$. The coincidence of the basis does not imply $\hat{f}(k) = \hat{g}(k)$: the coefficients depend on the expanded functions, not on the basis. Exactly the same way, the coincidence of the characters $\chi_k$ does not imply $G_k(A) = \lambda_k(A)$: the coefficients (orbital sums vs. first row) are different.

What is exactly true is this: The only case in which there is an exact relationship between $G_k$ the eigenvalues and the characteristic polynomial is through the characteristic polynomial, as stated in Proposition 11.15: the eigenvalues are the roots of the polynomial $\chi_A(\lambda) = \det(\lambda I - A)$, whose coefficients are the elementary symmetric functions $e_r(A) = \sum_{|I|=r} \det A[I,I] = \sum_{|I|=r} G_0^{(r)}(A[I,I])$. This relationship involves the fundamental modes of the vector minors of all orders $r = 1, \ldots, n$, not the crude modes $G_k(A)$ of total order $n$.

**Proposition 11.15 (Characteristic-Polynomial Bridge)**

The eigenvalues of $A$are not generally recovered from the crude modes $G_k(A)$, but from the characteristic polynomial

$$\chi_A(\lambda) = \det(\lambda I - A) = \lambda^n + \sum_{r=1}^{n} (-1)^r e_r(A)\, \lambda^{n-r},$$

where

$$e_r(A) = \sum_{|I|=r} \det\ A[I,I]$$

is the $r$-th elementary symmetric function of the eigenvalues, equal to the sum of all principal minors of order $r$.

Proof: The eigenvalues of $A$ are, by definition, the roots of

$\chi_A(\lambda) = \lambda^n - e_1(A)\lambda^{n-1} + \cdots + (-1)^n e_n(A),$

where $e_r(A)$ is the $r$-th elementary symmetric function of the eigenvalues. By the classical formula for the coefficients of the characteristic polynomial, $e_r(A) = \sum_{|I|=r} \det\ A[I,I]$, that is, $e_r(A)$ it is the sum of all principal minors of order $r$.

By Observation 5.2, the fundamental mode of each principal subblock satisfies $G_0^{(r)}(A[I,I]) = \det\ A[I,I]$.: Therefore:

$$e_r(A) = \sum_{|I|=r} G_0^{(r)}(A[I,I]).$$

Thus, the coefficients of $\chi_A$ are recovered from the hierarchy of fundamental modes of the vector minors of orders $r = 1, \ldots, n$. In contrast, each crude mode $G_k(A)$ is homogeneous of degree $n$ (Theorem 3.6), so it cannot by itself replace the collection $e_1(A), \ldots, e_n(A)$, whose degrees vary between 1 and $n$. The complete spectral information does not reside in the crude modes of the total order, but in the hierarchy of vector minors.

### 11.11. Final Classification: Intrinsic vs. Convention-Dependent

| **Category** | **Object** | **Condition** |
|---|---|---|
| **Absolute intrinsic** | $G_0(A) = \det(A)$,$R_\infty(A)$ | Without restriction |
| **Low intrinsic $\phi_{\text{can}}$** | $D_{\phi_{\text{can}}}(A)$, $\parallel D \parallel_2$, $V_2$, $H$, $P_A(q)$,$G(z;A)$ | Fixed$\sigma^*(1) = 1$ |
| **Global gauge covariant** | $D_\phi(A)$ | $D_{\phi'} = U_\beta P_u D_\phi$ |
| **Non-intrinsic low free orbital gauge** | Power distribution between modes, $V_2$,$H$ | Without fixing $\phi$ |

### 11.12. Limitations and Open Directions

Despite the structural development presented here, several fundamental questions remain open.

First, the injectivity of the map

$A \mapsto D_\varphi(A)$

is not yet characterized.

Second, the inverse synthesis problem, whether a prescribed modal vector can arise from a matrix, remains unresolved.

Third, the probabilistic behavior of modal energy and orbital entropy for random matrix ensembles has not yet been studied.

These directions suggest that the vector determinant should be viewed as the beginning of a broader research program rather than a completed classification theory.

**12. Conclusions**

The formalism developed in this work rigorously establishes a precise algebraic kernel for the ARE Method at its spectral level. This kernel consists of the space of phased monomials, the orbital decomposition induced by the action of $C_n$ on $S_n$, the orbital sums $\Lambda_r(A)$, their Fourier modes, $G_k(A)$ and the associated vector $D_\varphi(A)$. Within this framework, the classical determinant appears exactly as the fundamental mode:

$G_0(A) = \det(A)$.

The main conceptual contribution of the manuscript is to show that the classical determinant can be situated within a finite orbital-Fourier hierarchy:

$M_\varphi(A) \;\longrightarrow\; \{\Lambda_r(A)\} \;\Leftrightarrow\; D_\varphi(A) \;\longrightarrow\; G(z;A).$

Here, the map $M_\varphi(A) \to \{\Lambda_r(A)\}$ adds information by orbital classes and may lose internal structure; the map $\{\Lambda_r(A)\} \Leftrightarrow D_\varphi(A)$ is an exact Discrete Fourier Transform; and the map $D_\varphi(A) \to G(z;A)$ produces a finite trigonometric interpolation that extends the discrete modes to a continuous parameter.

From an axiomatic point of view, the central result must be interpreted precisely. The mode $G_0(A)$ coincides with the classical determinant and retains its alternating character. In contrast, the modes $G_k(A)$, $k \neq 0$, are multilinear observables associated with the orbital structure, but in general they are not alternating and, therefore, do not constitute determinants in the classical sense. Similarly, the values $G(\xi;A)$ for $\xi \in \mathbb{R}$ should be understood as finite trigonometric interpolations of the orbital modes, not as a canonical axiomatic extension of the classical determinant.

Among the established results are: the exact recovery of the determinant by means of the readout functional; Corollary 2.6 (Uniqueness of the Readout Functional); the multilinearity of all modes; the Hermitian symmetry for real matrices; the orbital and continuous Parseval identities; the vector Jacobi and Laplace formulas; the vector Hadamard inequality; and the degree-based obstruction, in general, to identifying the crude orbital modes with the eigenvalues of circulant matrices without introducing an additional nonlinear normalization.

The manuscript also honestly delimits several current limitations of the formalism. In particular, the theory of vector minors must always be developed in the symmetric group of the corresponding order; the dependence on the orbital gauge, outside the global affine case, does not reduce to a uniform covariance independent of the matrix; and the exact bridge between orbital structure and classical spectrum of circulants requires descending to the level of the characteristic polynomial and the minors of different orders.

Consequently, the work should not be presented as a closed theory, but rather as a rigorous framework already established at its core and open in several natural directions: classification of zeros of the orbital polynomial on the unit circle, inverse theory of spectral synthesis, fine probabilistic study of $V_2(A)$ low Gaussian models, complete formalization of vector minors, and possible extensions to tensors and hyperdeterminants. In this sense, the article provides a formal starting point for further research, but at the same time clearly identifies the questions that remain open.

**Appendix A**

**Proposition A.1 (Absence of uniform law: constructive verification)**

We complete the previous observation with an explicit calculation showing that $G_k(PA)$cannot, in general, be expressed as $\lambda_k(P)G_k(A)$independently of $A$.

Let $n = 3$, $\tau = (123)$ (cyclic rotation of rows), $\varepsilon(\tau) = +1$. The formula in Proposition 3.13 uses $\phi(\pi \circ \tau)$ (not $\pi \circ \tau^{-1}$). The compositions $\pi \circ \tau$ with $\tau = (123)$ for the permutations with $m_\pi(A) \neq 0$ are:

| $\pi$ | $\pi \circ \tau$ | $\phi(\pi \circ \tau)$ | $m_\pi(A)$ |
|---|---|---|---|
| $id$ | $(123)$ | $2$ | $+8$ |
| $(23)$ | $(13)$ | $2$ | $-2$ |
| $(12)$ | $(23)$ | $0$ | $-2$ |

For $\phi((13))$: the canonical representative of the orbit of $(13)$ satisfies $\sigma^*(1) = 1$; as $(13)(1) = 3 \neq 1$, we seek $r^* \equiv (13)^{-1}(1) - 1 = (3) - 1 = 2$ (mod3), that is $\sigma^* = (13) \circ \rho^2$. Now $(13) \circ \rho^2 : 1 \mapsto \rho^2(1) = 3 \mapsto (13)(3) = 1$. Then $\phi((13)) = 2$.

$$G_1(P_\tau A) = \varepsilon(\tau) \sum_{\pi} \varepsilon(\pi)\ m_\pi(A)\, \omega^{\phi(\pi\circ\tau)} = (+1)[8 \cdot \omega^2 + (-2) \cdot \omega^2 + (-2) \cdot \omega^0]$$

$$= 6\omega^2 - 2.$$

*With* $\omega^2 = -\frac{1}{2} - i\frac{\sqrt{3}}{2}$:

$$G_1(P_\tau A) = 6\left(-\frac{1}{2} - i\frac{\sqrt{3}}{2}\right) - 2 = -3 - 3i\sqrt{3} - 2 = -5 - 3\sqrt{3}\, i$$

This exactly matches the value of the Observation of Proposition 3.13. It is confirmed that $G_1(P_\tau A) = -5 - 3\sqrt{3}\,i \neq G_1(A) = 7 - i\sqrt{3}$ and $G_1(P_\tau A) \neq -G_1(A)$, verifying the absence of uniform law.